\documentclass[11 pt]{amsart}
\usepackage[all]{xy}
\usepackage{palatino}
\usepackage{mathtools}
\usepackage{latexsym, amsmath, amssymb, longtable, booktabs, amscd, microtype, booktabs, cases}
\usepackage[square]{natbib}
\usepackage{graphicx}
\usepackage[dvipsnames]{xcolor}
\usepackage{caption}
\usepackage{tcolorbox}
\usepackage{enumerate}
\usepackage{multirow}
\usepackage{tikz}
\usepackage{tikz-cd}
\usepackage{amsthm}

\theoremstyle{plain}
\newtheorem{theorem}{Theorem}[section]

\newtheorem{propn}{Proposition}[section]
\newtheorem{lem}[theorem]{Lemma}
\newtheorem{corollary}{Corollary}[section]

\newtheorem{mainthm}{Theorem}

\newtheorem{maincor}{Corollary}

\theoremstyle{definition}

\newtheorem{remark}{Remark}[section]

\usetikzlibrary{matrix,arrows,decorations.pathmorphing}
\usepackage{collectbox}
\usepackage{amsxtra}
\usepackage[colorinlistoftodos,prependcaption, textsize=tiny]{todonotes}
\reversemarginpar

\usepackage{enumitem}
\usepackage{mathabx}
\usepackage[square]{natbib}

\usepackage{color}
\usepackage{amsfonts}
\usepackage{geometry}
\usepackage{mathrsfs}
\usepackage{amsthm}
\usepackage{amsopn,xcolor}
\usepackage{amscd}
\usepackage{exscale, verbatim}
\usepackage[
  colorlinks=true,
  linkcolor={[RGB]{20,70,110}},
  citecolor=olive,
  urlcolor=black, linktoc=all]{hyperref}
\usepackage{bbold, bbm}
\bibpunct{[}{]}{,}{n}{}{;} 

\numberwithin{equation}{section}
\begin{document}
\renewcommand{\baselinestretch}{1.25}

\AtBeginEnvironment{theorem}{\setlist{before=\leavevmode, nosep}}

\newcommand{\abs}[1]{\left| #1 \right|}
\newcommand{\Q}{\mathbb{Q}}
\newcommand{\gln}{\mathrm{GL}_2(\mathbb{F}_\ell)}
\newcommand{\sln}{\mathrm{SL}_2(\mathbb{F}_\ell)}
\newcommand{\p}{\mathfrak{p}}
\newcommand{\Set}[1]{\{ #1 \}}
\newcommand{\gal}{\mathrm{Gal}(\overline{K}/K)}

\makeatletter
\def\imod#1{\allowbreak\mkern10mu({\operator@font mod}\,\,#1)}
\makeatother
\title{An explicit version of Serre's Open Image Theorem for Galois Number Fields}

\author{Clifford J Rodrigues}
\email{clifford.rodrigues@students.iiserpune.ac.in}
\address{Department of Mathematics, 
Indian Institute of Science Education and Research Pune,
Dr. Homi Bhabha Road, Pune 411008, INDIA.}

\author{Debargha Banerjee}
\email{debargha.banerjee@gmail.com}
\address{Department of Mathematics, 
Indian Institute of Science Education and Research Pune,
Dr. Homi Bhabha Road, Pune 411008, INDIA.}

\author{Krishna Kishore}
\email{kishorekrishna@iittp.ac.in}
\address{Department of Mathematics and Statistics, 
Indian Institute of Technology (IIT) Tirupati
Yerpedu – Venkatagiri Road, Yerpedu Post,
Tirupati District, Andhra Pradesh. India, 517619.
}
\author{Sunil Kumar Mallick}
\email{mallicksunil200114@gmail.com}
\address{Department of Mathematics and Statistics, 
Indian Institute of Technology (IIT) Tirupati
Yerpedu – Venkatagiri Road, Yerpedu Post,
Tirupati District, Andhra Pradesh. India, 517619.}

\thanks{The first author's research has been supported by the NBHM PhD Scholarship. The second author is supported by ANRF grant ANRF/ARGM/2025/000421/MTR.}

\subjclass[2020]{Primary: 11G05, Secondary: 14H52}
\keywords{Elliptic Curves, Galois Representations, Serre's Open Image Theorem.}

\begin{abstract}
 In 1971, Serre showed that the mod $p$ Galois representation associated to a non-CM elliptic curve over a number field is surjective for all but finitely many primes $p$. That is, there exists a constant $c(E,K)$, for $E$ a non-CM elliptic curve over a number field $K$, such that $\rho_{E,p}:\gal\to \mathrm{Aut}(E[p])$ is surjective for all primes $p>c(E,K)$. We give a bound for this constant using an effective version of the Faltings isogeny theorem for specific number fields. We show that, for $K=\mathbb{Q}(\sqrt{D})$ with $|D|<50$ and $K$ a formally real number field, given a non-CM elliptic curve $E$ over $K$ we have $c(E,K)\le 2C_1(K)\log N_{K/\mathbb Q}\operatorname{rad}\left( 2\mathcal O_K\mathfrak N_{E}\right )+2 C_2(K)$, where $\mathfrak N_E$ denotes the conductor ideal of $E$, and the constants $C_1(K)$ and $C_2(K)$ can be described explicitly. We also show that for primes greater than this bound, the image cannot be contained in the normalizer of a non-split Cartan subgroup. This works for any Galois number field. Our proof assumes the Generalized Riemann Hypothesis (GRH).
\end{abstract}
\maketitle

\section{Introduction}
Elliptic curves are among the important objects of study in number theory. To understand them, we study various associated quantities, such as their groups of rational points, the primes of good reduction, and their associated Galois representations, to name a few. The last quantity, in particular, provides a lot of information about the original curve. For example, the Neron-Ogg-Shafarevich theorem states that the primes at which the Galois representation is unramified correspond to the primes of good reduction for the elliptic curve. Understanding the image of this representation is therefore an equally useful objective. One of the strongest theorems along these lines is Serre's open image theorem, which we now describe.\\
Let $E$ be an elliptic curve without complex multiplication over a number field $K$. Let $G_K:=\mathrm{Gal}(\overline{K}/K)$ denote the absolute Galois group of $K$. A celebrated theorem of Serre (\cite[Th\'eorem\`e 2, p. 293]{Serre1971}) states that the associated Galois representation
\begin{equation*}
    \rho_{E,\ell}:G_K\rightarrow \text{Aut}(E[\ell])\simeq \gln
\end{equation*}
is surjective for all but finitely many primes $\ell.$ Therefore, there exists a smallest positive integer $c(E)$, depending on $E$, such that for all primes $\ell>c(E)$, the map $\rho_{E,\ell}$ is surjective. Serre managed to prove that $c(E)=11$ would work for semi-stable elliptic curves over $\Q$ \cite[Corollaire 2.,p 309]{Serre1971}.\\
In 1993, Masser and Wustholtz \cite{masser} gave a bound $c(E)<c_2(\mathrm{max}(d,h))^\gamma$ where $d$ denotes the degree of the number field $K$, $h$ denotes the Weil height of the j-invariant of $E$ and $c,\gamma$ are absolute constants. 
One can also look at the representation on the Tate module $T_\ell(E)$ of $E$, that is the map $\overline{\rho}_{E,\ell}:G_K\rightarrow \text{Aut}(T_\ell(E))\simeq \mathrm{GL}_2(\mathbb Z_\ell)$. We may further collect all the $\overline{\rho}_{E,\ell}$ into one representation, namely, $\rho_\infty:G_K\to \mathrm{GL_2}(\hat{\mathbb Z}_\ell)$. An equivalent form of Serre's theorem says that the index of $\rho_\infty(G_K)$ in $\mathrm{GL_2}(\hat{\mathbb Z}_\ell)$ is bounded by a constant that depends only on $E$ and $K$. Towards this, Davide Lombardo \cite{lombardo} gives an unconditional bound on the index of the adelic image of the Galois representation.\\
Assuming the Generalized Riemann hypothesis (GRH), allows one to give an explicit version of the Chebotarev density theorem which, in turn, helps us give better bounds for $c(E)$. For example, Serre showed that for elliptic curves over the rationals \cite[Th\'eorem\`e 22, p. 196]{serre2} one has 
\begin{equation*}
    c(E)\leq C_4(\log N_E)(\log\log2N_E)^3
\end{equation*}
In 2014, Larson and Vaintrob \cite{larsonvain} show that $c(E)\leq C_3 \log N_E$, although the constant $C_3$ is not effectively computable.\\
In 2023, Mayle and Wang \cite{maylewang} give an upper bound for $c(E)$ for non-CM ellpitic curves over $\Q.$ They use an effective version of Chebotarev's density theorem by Bach and Sorensen \cite{bach1996explicit}. This version of the explicit Chebotarev also holds for number fields. Using a better quotient, Chen and Swidinksy\cite{Chen_Swidinsky_2026} have been able to improve the bound given by Mayle and Wang.\\
 For number fields, the situation remains largely unknown. For example, Mazur's cyclic isogeny theorem \cite[Theorem 1]{mazur} completely classifies the primes for which an elliptic curve can have a cyclic isogeny of that degree. This in turn allows one to rule out the cases for which $\rho_{E,\Q}(G_\Q)$ is contained in a Borel subgroup (see Theorem \ref{borel iso}). However, for quadratic fields (see Banwait \cite{banwait}), the cyclic isogenies have been determined. In fact, an algorithm is provided that takes a quadratic field as input and outputs the cyclic isogeny primes. For higher-degree fields, the situation is unknown.
 \subsection{Main Results} We extend the results of Mayle and Wang to Galois number fields. The main challenge with number fields is to keep track of the various prime ideals lying over a given rational prime. This leads to subtle differences in the statements of the theorems compared with Mayle and Wang. Furthermore, many of the constants associated with elliptic curves now turn out to be ideals. For example, the conductor of an elliptic curve is replaced by the conductor ideal. This still tracks the prime ideals of bad reduction.\\
 Our main focus will be to deal with the case in which the image is contained in the normalizer of a non-split Cartan subgroup but not fully contained in said non-split Cartan subgroup (see Theorem \ref{subgroup classification}, case (v).) 
 \begin{mainthm}[Theorem \ref{5.4mainthm}]\label{main thm simple}
     Let $E$ be a non-CM elliptic curve, with conductor ideal $\mathfrak N_E$, over a Galois number field $K$. For any prime number $\ell$ such that,
     \begin{equation*}
         \ell>2C_1(K)\log N_{K/\mathbb Q}\operatorname{rad}\left( 2\mathcal O_K\mathfrak N_{E}\right )+2 C_2(K) ,
     \end{equation*}
     the image of the Galois representation $\rho_{E,\ell}$ is not contained in the normalizer of a non-split Cartan subgroup, where
\begin{align*}
    C_1(K)=& 510.352+513.024n_{K},\\
    C_2(K)=& 513.024\bigl(1+\frac{n_{K}}{4}\log d_K
+n_{K}\log(384n_{K})\bigr)
+1324.8n_{K}+45.7,
\end{align*}
and $d_K$ denotes the absolute discriminant of the number field $K.$
 \end{mainthm}

The proof of this theorem allows us to give a bound for $c(E,K)$, where $K=\Q(\sqrt{D})$ is a totally real number field with $|D|<50$. We assume that $K$ is quadratic to apply Banwait's classification of cyclic isogenies \cite{banwait}. This, in turn, allows us to bound the prime numbers for which the image is contained in a Borel subgroup. This is also where we get the restriction of $|D|<50$. Similarly, we assume that $K$ is totally real to restrict the primes for which the image lies in a non-split Cartan subgroup. Note that the classification of maximal subgroups (Theorem \ref{subgroup classification}) makes a distinction between being contained in a non-split Cartan and being contained in the normalizer of the non-split Cartan (see Theorem \ref{subgroup classification} (iii) and (v)) We assume GRH to provide an effective version of Faltings Isogeny Theorem (which in turn uses an effective version of the Chebotarev density theorem from \cite{bach1996explicit}, which is where the assumption is needed) Putting it all together, we have the following:

    \begin{maincor}\label{total theorem}
    Let $K=\Q(\sqrt{D})$ be a formally real quadratic number field with $|D|<50$ and let $E$ be an elliptic curve over $K$ without complex multiplication. Let $\rho_{E,\ell}:G_K\rightarrow \gln$ be the associated Galois representation. Assume GRH and suppose the image of $\rho_{E,\ell}$ is not contained in the normalizer of a split Cartan subgroup of $\gln$. Let $c(E,K)$ denote the largest prime for which the image of the Galois representation fails to be surjective. Then
    \begin{equation*}
        c(E,K)\le 2C_1(K) \log N_{K/\mathbb Q}\operatorname{rad}\left( 2\mathcal O_K\mathfrak N_{E}\right )+2 C_2(K) ,
    \end{equation*}
    where the constants $C_1(K)$ and $C_2(K)$ are the same as in Theorem \ref{main thm simple}.
\end{maincor}
 For a prime ideal $\mathfrak p$ of $K$, we denote by $a_\mathfrak p(E)$ the constant
\begin{equation*}
    a_\mathfrak p(E)=N(\mathfrak p)+1-\#E(\mathbb F(\mathfrak p)),
\end{equation*}
where $N(\mathfrak p)$ denotes the norm of the ideal $\mathfrak p$, $\mathbb  F(\mathfrak p)$ denotes the residue field at $\mathfrak p$ and $\#E(\mathbb F(\mathfrak p))$ denotes the number of point of the elliptic curve $E$ in the reduction modulo $\mathfrak p.$\\
To prove Theorem \ref{main thm simple}, we prove an effective version of Faltings isogeny theorem over number fields. This theorem is central to all that follows. To adapt this result from Mayle and Wang \cite[Theorem 2, p. 1401]{maylewang} to number fields, we first prove a series of number theoretic lemmas in $\S 4.2$. This allows us to construct an extension in which a given collection of prime ideals ramifies. We can use an extension to find a prime ideal of good reduction by ensuring that all the bad-reduction primes ramify. Combining this with the effective Chebotarev density theorem by Bach and Sorenson \cite{bach1996explicit} and the generalized discriminant bound (Lemma \ref{bound disc}), we are able to prove the following:

  \begin{mainthm}[Theorem \ref{Eff iso for K}]
Assume GRH. Let $K$ be a Galois number field, and let $E_1,E_2/K$ be elliptic curves without complex multiplication. Suppose that $E_1$ and $E_2$ are not $K$-isogenous. Then there exists a prime ideal $\mathfrak p\subseteq \mathcal O_K$ of good reduction for both $E_1$ and $E_2$ such that
$$
a_{\mathfrak p}(E_1)\neq a_{\mathfrak p}(E_2),
$$
and satisfying
$$
N_{K/\mathbb Q}\mathfrak p
\le
\left(
C_1(K)
\log N_{K/\mathbb Q}
\operatorname{rad}
\left(
2\mathcal O_K\mathfrak N_{E_1}\mathfrak N_{E_2}
\right)
+
C_2(K)
\right)^2,
$$
where $\mathfrak N_{E_i}$ is the conductor ideal of $E_i$, and $C_1(K),C_2(K)$ are constants depending only on $K$.
\end{mainthm}
    \subsection{Related Work}
    Le Fourn \cite{lefourn} proves a uniform surjectivity result for non-CM $\Q$-curves that are defined over an imaginary quadratic field. \\\
    Lemos \cite{lemos} shows the existence of a constant $C_{K,d}$ for which the surjectivity of the representation $\mathbb P\rho_{E,p}$ is established in various cases, where $E$ is a $\Q$-curve defined over a quadratic field $K$. The constant is not made explicit, except in \cite[Theorem $1.10$]{lemos}.
   \subsection{Outline of the Paper}
   The outline of the proof that we follow is due to Mayle and Wang \cite{maylewang}.\\
   In section 2, we introduce notations and describe the classification of maximal subgroups for $\gln$. Using this classification, for each maximal subgroup $H$, we describe conditions for a prime number $\ell$ under which the image of the galois representation $G_\ell$ is not contained in $H$.\\
   Our main focus of this paper is to give a condition on $\ell$ for which $G_\ell$ is not contained in the normalizer of a non-split Cartan subgroup. Towards this end, Section 3 provides bounds for the norm of a prime ideal $\mathfrak p$. This step uses an effective version of the Chebotarev density theorem established by Bach and Sorenson \cite{bach1996explicit}. It is here that we need to assume the Generalized Riemann hypothesis.\\
   Section 4 establishes an effective version of Faltings Isogeny Theorem for number fields. That is, given two elliptic curves $E$ and $E'$ over a number field $K$ that are not $K$-isogenous, we show that there exists a prime ideal $\mathfrak p$ for which $a_p(E)\neq a_p(E')$. In fact, by using the bounds established in Section 3, we are able to provide a bound on the norm of the prime ideal $\mathfrak p$ thus obtained.\\
   Finally, section 5 completes the proof of the main theorem. We first assume that the image $G_\ell$, of an elliptic curve $E$ over $K$, is contained in the normalizer of the non-split Cartan subgroup but not completely contained in the non-split Cartan subgroup. From this, we construct a quadratic character $\chi_D$. We twist the elliptic curve $E$ with respect to this character to get another elliptic curve $E_D$. After this, we show that these elliptic curves are non-isogenous over $K$ and hence, by the effective version of Faltings isogeny theorem established in Section 4, we obtain a prime ideal $\mathfrak p$ for which $a_\mathfrak p(E)\neq a_\mathfrak p(E_D)$ and whose norm satisfies a certain bound. Finally, by applying the Hasse bound, we obtain the desired bound for the prime $\ell$.
   \subsection{AI Disclosure}
   In preparing this paper, ChatGPT and DeepSeek were used to proofread and to help with finding references for some results. The authors take full responsibility for the contents presented within this paper.
    \section{Preliminaries}
    In this section, we describe the conditions on the field $K$ and on the prime number $\ell$ for which the image fails to be surjective. For this, we appeal to the classification of subgroups of $\gln$, which dates back to Dickson \cite{dickson} (See Theorem \ref{subgroup classification}) Using this, we do a case-by-case analysis by assuming the image is contained in the various subgroups and obtaining a bound on the prime $\ell$. Combining all these results together with Theorem \ref{main thm simple}, we get Corollary \ref{total theorem}.
    \subsection{Notations}
    We use the same notations as in Mayle and Wang \cite[\S 2, p.1402]{maylewang}. Let $\ell$ be an odd prime. For $\varepsilon \in \mathbb{F}_\ell$, let $\left( \frac{\varepsilon}{l} \right)$ denote the Legendre symbol. Let
    \begin{align*}
        C_s(\ell):&=\left\{\begin{pmatrix}
            a&0\\
            0&b
        \end{pmatrix}:a,b\in \mathbb{F}_\ell,\;ab\neq 0\right\}\quad \text{and} \\  C_{ns}(\ell):&=\left\{\begin{pmatrix}
            a&\varepsilon b\\
            b&a
        \end{pmatrix}\in\mathrm{GL}_2(\mathbb F_\ell):a,b,\varepsilon\in \mathbb{F}_\ell,\;\left(\frac{\varepsilon}{\ell}\right)\neq 1\right\},
    \end{align*}
    denote the standard split and non-split Cartan subgroups, respectively. By $\left( \frac{\varepsilon}{\ell} \right)\neq 1$, we mean that $\varepsilon$ is a quadratic non-residue in $\mathbb F_\ell^\times.$ 
    Let $ C_s^+(\ell)$ and  $C_{ns}^+(\ell)$ denote the normalizers in $\gln$ of the split Cartan and non-split Cartan subgroups, respectively. Let
    \begin{equation*}
        B(\ell):=\left\{\begin{pmatrix}
            a&b\\
            0&d
        \end{pmatrix}:a,b,d\in \mathbb{F}_\ell,\;ad\neq 0\right\}, 
        \end{equation*}
    denote the standard Borel subgroup in $\gln$. Let $
    \overline{G}$ denote the image of a subgroup $G$ of  $\mathrm{GL}_2(\mathbb F_\ell)$ under the natural surjection $q:\mathrm{GL}_2(\mathbb F_\ell)\to \mathrm{PGL}_2(\mathbb F_\ell)$. Let $A_n$ (resp. $S_n$) denote the alternating (resp. symmetric) groups in $n$ letters. We then have the following classification due to Dickson:
     \begin{theorem}[Dickson]\label{subgroup classification}
        Let $\ell$ be an odd prime. If $G\subseteq \gln $ is a subgroup, then either
        \begin{enumerate}[label={(\roman*)},itemindent=1em]
            \item $G$ contains $\sln ,$
            \item $G$ is conjugate to a subgroup of $B(\ell)$,
            \item $G$ is conjugate to a subgroup of $C_{ns}(\ell)$,
            \item $G$ is conjugate to a subgroup of $C_{s}^+(\ell)$ but not to any subgroup of $C_s(\ell)$,
            \item $G$ is conjugate to a subgroup of $C_{ns}^+(\ell)$ but not to any subgroup of $C_{ns}(\ell)$, or
            \item $\overline{G}$ is isomorphic to $A_4,S_4$ or $A_5.$
        \end{enumerate}
    \end{theorem}
    \begin{proof}
       The reader may refer to section 2 of Serre's original paper \cite[\S  2]{Serre1971} for a proof of this result.
    \end{proof}
   Let $G_\ell$ be the image of the Galois representation $\rho_{E,\ell}:G_K\rightarrow \gln$. Our strategy is to describe constraints on primes for which the image $G_\ell$ is contained in the various maximal subgroups appearing in Theorem \ref{subgroup classification}. That is, we provide a bound $N$ such that, for primes $\ell>N$, the subgroup $G_\ell$ will not be contained in any maximal subgroup. Taking the maximum of each such bound will give us an estimate for $c(E)$ in Serre's Open Image Theorem.

    \subsection{$G_\ell$ contains $\sln$}
   Suppose that the image of the Galois representation contains $\sln.$ We have the following: 
    \begin{propn}
        Suppose that $\ell$ is not ramified in $K$. If $G_\ell$ contains $\sln$, then $G_\ell$ is $\gln$, i.e. $\rho_{E,\ell}$ is surjective.
    \end{propn}
   
    \begin{proof}
         Using the Weil pairing on $E$, the composition 
        \begin{equation*}
            \det\circ \rho_{E,\ell}:G_K\rightarrow \mathbb{F}_\ell^\times
        \end{equation*}
        is the mod $\ell$ cyclotomic character, by a similar calculation as in the proof of \cite[Proposition 8.6, p. 99]{Silverman2009}, and is surjective when $\ell$ is unramified in $K/\mathbb Q$. 
Suppose that $\sln\subset G_\ell$. The image of $\sln$ under the determinant map is $1$, but by surjectivity of the cyclotomic character we see that $\sln\subsetneq G_\ell$. Therefore, by Theorem \ref{subgroup classification}, we will have $G_\ell=\gln$.
    \end{proof}
    \subsection{$G_\ell$ is conjugate to a subgroup of  $B(\ell)$}
      In this subsection, we give a condition under which the image $G_\ell$ is contained in $B(\ell)$. We begin by relating $G_\ell$ to the existence of isogenies of the associated elliptic curve. This theorem is well-known (see, for example, \cite[\S 1, p.314]{nuriavila}) but due to the lack of a proper reference, we provide a proof
    \begin{propn}\label{borel iso}
        $G_\ell$ is contained in a $B(\ell)$ if and only if the elliptic curve $E$ admits an isogeny, defined over $K$, of degree $\ell$.
    \end{propn}
    \begin{proof}
        Suppose that $\rho_{E,\ell}(G_K)$ is contained in a Borel subgroup. Without loss of generality, we may assume it is contained in $B(\ell)$. Thus $E[\ell]$ has a one dimensional invariant subspace of order $\ell$, say $C.$ Using this subgroup, we can construct a quotient isogeny
        \begin{equation*}
            \phi: E\rightarrow E/C
        \end{equation*}
        defined over $K$,
        of degree $\deg \phi= \# \ker \phi= \# C=\ell.$
      Hence $E$ admits an isogeny of degree $\ell.$\\
    Conversely, assume that $\varphi:E\rightarrow E'$ is an isogeny of degree $\ell.$ The order of the subgroup $\ker\varphi$ equals the degree of $\varphi$, which is $\ell.$ 
    Hence
    \begin{equation*}
        \ker\varphi\subset E[\ell].    \end{equation*}
        Since $\varphi$ is defined over $K$, we have $\rho_{E,\ell}(G_K)(\ker\varphi)\subset \ker\varphi.$
        Choosing a generator of the group $\ker\varphi$ as one of the basis elements for $E[\ell]$, we see that 
        \begin{equation*}
            \rho_{E,\ell}(G_K)\subset \begin{pmatrix}
                *&*\\
                0&*
            \end{pmatrix}
        \end{equation*}
        Thus, the image is contained in a Borel subgroup.
        \end{proof}
       Thus, if we know the primes for which the elliptic curve admits a rational isogeny of that degree, we will have found the primes for which the image of the associated Galois representation is contained in a Borel subgroup. According to Theorem 1.9 of Banwait \cite{banwait}, for primes greater than the integer $163$, should ensure that the image is not contained in a Borel subgroup for $\Q(\sqrt{D})$ for $|D|<50.$ 

        \subsection{$G_\ell$ is contained in an Exceptional subgroup}
        Let $K$ be a number field, $\ell$ a rational prime number and $\nu$ a prime of $\mathcal{O}_K$ lying over $\ell$. Let $E$ be an elliptic curve over $K$. From Serre's analysis of the tame inertia subgroups (refer to sections 1.11 and 1.12 in Serre \cite[ Proposition 11, Proposition 12]{Serre1971}), we have the following:
        \begin{enumerate}
            \item If $E$ has ordinary good reduction at $\nu$ or multiplicative bad reduction, then $\rho_{E,\ell}(I_w)$ has order $\ell-1$ or $\ell(\ell-1)$ and one can represent it in matrices of the form
            \begin{equation*}
                \begin{bmatrix}
                    *&0\\0&1
                \end{bmatrix} \quad \text{or}\quad
                \begin{bmatrix}
                    *&*\\0&1
                \end{bmatrix}.
            \end{equation*}
            \item If $E$ has supersingular good reduction at $\nu$, then $\rho_{E,\ell}(I_w)$ is cyclic of order $\ell^2-1.$
         \end{enumerate}
         Let $\tilde{G_\ell}$ be the projection of $G_\ell=\rho_{E,\ell}(G_K)$ in the group $\mathrm{P}\gln$. In case 1, $G_\ell$ contains a split semi-Cartan subgroup whose projection in $\mathrm{P}\gln$ is a cyclic subgroup of order $\ell-1$ (Discussion of part (a)  in \cite[\S 2.1]{Serre1971}). In case ($2$), $G_\ell$ contains a nonsplit Cartan subgroup whose projection in $\mathrm{P}\gln$ is a cyclic subgroup of order $\ell+1$.\\
         Since the subgroups $A_4,S_4$ and $A_5$ have elements of order at most 6, choosing a prime number $\ell$ greater than $7$ implies that $\tilde{G_\ell}$ will not be equal to these exceptional subgroups.
        \subsection{$G_\ell$ is conjugate to a subgroup of $C_{ns}(\ell)$}
    We finally consider the case where the image lies in a non-split Cartan subgroup. It is here that we crucially require the field to be a formally real number field.
    \begin{lem}
        If $G_K$ contains an involution, then the image cannot be contained in a non-split Cartan subgroup.
    \end{lem}
    \begin{proof}
        Let $K$ be a field such that $\gal$ contains an involution (for example, $K$ is a formally real number field \cite[ Kapitel I, \S 6, Korollar 2]{manfredclaus}.)\\
        Then $\rho_\ell(G)$ contains an element with eigenvalues $\pm 1$. However, $C_{ns}(\ell)$ does not contain any element with these eigenvalues. To see this, we note that a general element of $C_{ns}(\ell)$ is of the form
        \begin{equation*}
            A=\begin{pmatrix}
                a&\varepsilon b\\
                b& a
            \end{pmatrix}
            \quad a,b\in \mathbb{F}_\ell\; \text{and $\varepsilon$ a nonsquare in $\mathbb{F}_\ell$.}
        \end{equation*}
        For $A$ to have eigenvalues $\pm 1$, we must have
        \begin{equation*}
            tr(A)=0\quad\text{and}\quad det(A)=-1. 
        \end{equation*}
        This would imply $2a=0$ and $a^2-\varepsilon b^2=-1$. Assuming $\ell\neq 2$, the first of these relations gives us $a=0$, which implies that we have $\varepsilon b^2=1,$ or $\varepsilon= (b^{-1})^2$. This is a contradiction to the fact that $\varepsilon$ was assumed to be a nonsquare in $\mathbb{F}_\ell.$ 
    \end{proof}
    To summarize, we have the following. 
     From Theorem \ref{subgroup classification}, all that is left to be classified is the case when
     \begin{enumerate}
         \item $G_\ell$ is conjugate to a subgroup of the normalizer of a split Cartan subgroup, but not conjugate to any subgroup of the split Cartan.
         
         \item $G_\ell$ is conjugate to a subgroup of the normalizer of a non-split Cartan subgroup, but not conjugate to any subgroup of the non-split Cartan.
     \end{enumerate}
     For rational numbers, case (1) has been shown to not be possible for $\ell\geq11$ and $\ell\neq 13$ by Bilu, Parent and Reboleddo \cite[Corollary 1.2, p.960]{biluparent}. In the paper, they show that $X_{\text{split}}(\ell)$ has no non-trivial non-CM points for $\ell\geq13$ and $\ell\neq 13$. The case $\ell=13$ was handled by Balakrishana, Dogra, Netan, M\"uller and Vonk \cite[Theorem 1.2, Corollary 1.3, p.888]{dogravonk} using quadratic Chabauty methods. To answer the question for number fields, one has to determine the non-CM $K$-rational points on $X_{\text{split}}(\ell)$. This is still an open problem. We do not deal with this case in this paper and only discuss case (2).\\
     When we do not wish to make a distinction between the split and non-split Cartan subgroups, we will just refer to it as a Cartan subgroup. For now we have the following theorem,
    \begin{theorem}\label{bound theorem}
        Let $E$ be an elliptic curve without complex multiplication over a formally real number field $K=\Q(\sqrt{D})$ with $0<D<50$. If $\ell$ is a prime number such that $\ell>163$,  then either $\overline{\rho}_{E,\ell}$ is surjective or $G_E(\ell)$ is conjugate to a subgroup of the normalizer of a Cartan subgroup but is not conjugate to any subgroup of the Cartan subgroup.
    \end{theorem}

\section{Image contained in normalizer of non-split Cartan but not contained in non-split Cartan}
Our goal will be to obtain a bound for which the image of the Galois representation can be contained in the normalizer of a non-split Cartan subgroup of $\gln$. The goal of this section will be to modify the effective Chebotarev theorem established in \cite{bach1996explicit} to get a version that allows one to avoid a prescribed set of prime ideals.

\subsection{Forms of the Chebotarev density theorem}

Let $K$ be a number field of degree $n_K$ over the rationals $\Q$, and $L$ be a finite normal extension  of $K$ such that $n_L = [L:\Q]$ and $n = n_L/n_K = [L:K]$, the degree of $L$ over $K$. Let $G = \mathrm{Gal}(L/K)$ be the Galois group of $L$ over $K$.
Let $d_K, $ resp. $d_L$ be the absolute discriminant of $K$, resp.  of $L$. 

Let $\mathfrak{p}$ be a prime ideal of $K$ unramified in $L$, and let the Artin symbol $\left( \frac{\p}{L/K} \right)$ (for more on the Artin symbol, refer to \cite[p. 106-109]{davidcox}) denote the conjugacy class of Frobenius automorphisms corresponding to prime ideals $\mathfrak P \mid \p$, where $\mathfrak P$ is a prime ideal of $L$. For each conjugacy class $C$ of $G$, we define
$$
\pi_C(x, L/K) = \abs{ \left\{\p : \p \textrm{ unramified in } L, \left( \frac{\p}{L/K} \right) = C,\; N_{K/\Q}(\p) \leq x \right\} }
$$
Let $\textrm{Li}(x)$ be the familiar logarithmic integral: 
$$
\textrm{Li}(x) = \int_2^x \frac{dt}{ \log t} \sim \frac{x}{\log x} ,\textrm{ as } x \to \infty
$$
By Generalized Riemann Hypothesis (GRH) we mean the assertion that the Dedekind zeta function of a number field has no zeros with real part greater than $1/2$. This hypothesis enables one to obtain the following strengthened version of a consequence of the Chebotarev density theorem \cite{lagarias}:

\begin{theorem}
    Assuming GRH, there exists an absolute constant $c > 0$ such that 
    $$
    \abs{ \pi_C(x) - \frac{\abs{C}}{\abs{G}} Li(x) } \leq c \frac{\abs{C}}{\abs{G}} \sqrt{x} (\log d_L + n_L \log x) 
    $$
    for all $x \geq 2$.
\end{theorem}
We use the following explicit version due to Bach and Sorenson \cite[Theorem 5.1, p. 1729]{bach1996explicit}
\begin{theorem}\label{bachthm}
    Assume GRH. Let $L/K$ be a Galois extension of number fields, with $L \neq \mathbb{Q}$. Let $d_L$ denote the absolute value of $L$'s discriminant. Let $n$ denote the degree of $L$ over $\mathbb Q$. Let $\sigma \in G$, the Galois group of $L/K$.
Then there is a prime ideal $\mathfrak{p}$ of $K$ with $\left(\frac{\mathfrak{p}}{L/K}\right) = \sigma$, of residue degree $1$, satisfying
\begin{equation*}
N\mathfrak{p} \leq (4 \log d_L + 2.5n + 5)^2.
\end{equation*}
\end{theorem}

\subsection{A modified version of the Effective Chebotarev density theorem}
We wish to modify the conclusion of Theorem \ref{bachthm}, so that we are able to choose a prime ideal $\mathfrak p$ that lies outside a prescribed set $S$. If we suppose that $S$ contains the prime ideals of bad reduction for an elliptic curve $E$, then the modified theorem would give us a prime ideal $\mathfrak p$ of good reduction that satisfies the conclusion of Theorem \ref{bachthm}. To do so, we first prove
\begin{lem}\label{prime ramify extn}
    Let $K$ be a number field and let
$S=\{\mathfrak q_1,\ldots,\mathfrak q_e\}$
be a finite set of nonzero prime ideals of $\mathcal O_K$. Then there exists an element $\alpha\in K^\times$ such that
$
v_{\mathfrak q_i}(\alpha)\equiv 1 \pmod 2
$
for $1\leq i\leq e$. For every $\alpha$ satisfying these valuation conditions, every prime ideal of $S$ ramifies in the extension $M=K(\sqrt{\alpha})$.
\end{lem}
\begin{proof}
    We first show the existence of such an $\alpha$ satisfying the prescribed valuation conditions. For each $1\leq i\leq e$, choose an $a_{i}\in \mathfrak {q}_{i}\setminus \mathfrak{q}_{i} ^{2}$. By the Chinese remainder theorem, there exists an $\alpha \in \mathcal{O}_{K}$ such that $\alpha\equiv a_{i}\pmod{\mathfrak{q}_{i}^{2}}$ for every $1\leq i\leq e$. 
    For each $i$, the inclusion $\alpha-a_{i}\in \mathfrak{q}_{i}^{2}$ gives $\alpha\in \mathfrak{q}_{i}.$ Moreover, if $\alpha\in \mathfrak{q}_{i}^{2},$ then $a_{i}=\alpha-(\alpha-a_{i})\in \mathfrak{q}_{i}^{2}$, contrary to the choice of $a_{i}.$ Consequently, $\alpha\in  \mathfrak{q}_{i}\setminus  \mathfrak{q}_{i}^{2}$ for every $1\leq i\leq e$. In particular, $\alpha$ is non-zero.
    
    We now prove that each prime ideal in $S$ ramifies in $M=K(\sqrt{\alpha})$. Fix an index $i,$ and suppose $(\alpha)=\mathfrak{q}_{i}^{r} \mathfrak{p}_{1}^{r_{1}}\cdots \mathfrak{p}_{t}^{r_{t}},$ where $r\ge 0$ and the $\mathfrak{p}_{j}$ are distinct nonzero prime ideals different from $\mathfrak{q}_{i}$. Since $\alpha\in \mathfrak{q}_{i}\setminus \mathfrak{q}_{i}^{2}, $ by \cite[Corollary 3.13, p. 51]{milne2020algebraic} $
v_{\mathfrak q_i}(\alpha)\equiv 1 \pmod 2
$
for every $1\leq i\leq e$. Now let $\alpha\in K^{\times}$ be any element such that $v_{\mathfrak q_i}(\alpha)$ is odd, then the prime ideal $\mathfrak q_i$ ramifies in $M/K.$ Indeed, fix $\mathfrak{q}_{i}\in S$ and put $r:=v_{\mathfrak{q}_{i}}(\alpha).$ Since $r$ is odd, $\alpha$ is not a square in $K.$ Consequently, $M=K(\sqrt{\alpha})$ is a quadratic extension of $K.$ Let $\mathfrak{Q}$ be a prime ideal of $\mathcal{O}_{M}$ lying above $\mathfrak{q}_{i}.$ Express $\mathfrak{q}_{i}\mathcal{O}_{M}=\mathfrak{Q}^{e}\mathfrak{Q}_{2}^{e_{2}}\cdots \mathfrak{Q}_{g}^{e_{g}},$ where $\mathfrak{Q}$ is distinct from $\mathfrak{Q}_{j},\ 2\le j\le g$ and we have $\alpha \mathcal{O}_{k}=\mathfrak{q}_{i}^{r} \mathfrak{p}_{1}^{r_{1}}\cdots \mathfrak{p}_{t}^{r_{t}},$ where $\mathfrak{p}_{j}$ are distinct from $\mathfrak{q}_{i}.$ Now extending the ideals to $\mathcal{O}_{M},$ we obtain $$ \alpha\mathcal{O}_{M}=(\mathfrak{q}_{i}\mathcal{O}_{M})^{r} (\mathfrak{p}_{1}\mathcal{O}_{M})^{r_{1}} \cdots (\mathfrak{p}_{t}\mathcal{O}_{M})^{r_{t}}.$$ So $\mathfrak{Q}$ occurs with exponent $er$ in $(\mathfrak{q}_{i}\mathfrak{O}_{M})^{r}.$ Moreover, $\mathfrak{Q}$ does not occur in any $\mathfrak{p}_{j}\mathcal{O}_{M}.$ Indeed, \cite[Lemma 3.33, p. 59]{milne2020algebraic} says that $\mathfrak{Q}\mid \mathfrak{p}_{j}\mathcal{O}_{M}$ if and only if $\mathfrak{Q}\cap \mathcal{O}_{K}=\mathfrak{p}_{j}.$ But $\mathfrak{Q}$ lies above $\mathfrak{q}_{i},$ so $\mathfrak{Q}\cap \mathcal{O}_{K}=\mathfrak{q}_{i},$ which is not equal to $p_{j}.$ Hence $v_{\mathfrak{Q}}(\alpha)=er$ in $\alpha \mathcal{O}_{M}.$ Since, $\alpha \mathcal{O}_{M}=(\sqrt{\alpha}\mathcal{O}_{M})^{2},$ so the exponent of $\mathfrak{Q}$ in $\alpha \mathcal{O}_{M}$ is even, so $e$ must be even. Since $M/K$ is Galois, so from \cite[Theorem 3.34, p. 59]{milne2020algebraic}, we conclude that $e\mid2,$ which implies that $e=2>1.$ So $\mathfrak{q}_{i}$ ramifies in $M/K.$
\end{proof}
Using this Lemma, we prove a modified version of the Bach and Sorenson theorem, namely
\begin{lem}(Bound for the norm of a prime ideal)\label{bach norm prime}
Assume GRH. Let $K$ be a number field and let $L/K$ be a finite Galois extension. Let
$C\subseteq \mathrm{Gal}(L/K)$ be a nonempty subset which is closed under conjugation.
Let
$S=\{\mathfrak q_1,\ldots,\mathfrak q_e\}$
be a finite set of nonzero prime ideals of $\mathcal O_K$. Fix any $\alpha\in K^\times$ satisfying
$
v_{\mathfrak q_i}(\alpha)\equiv 1 \pmod 2
$
for every $1\leq i\leq e$ (whose existence follows from Lemma \ref{prime ramify extn}) and put 
$
M:=K(\sqrt{\alpha})
$
and
$
\widetilde L:=LM.
$
Then there exists a  prime ideal $\mathfrak p\subseteq \mathcal O_K$ such that
$
\mathfrak p\notin S,
$
$\mathfrak p$ is unramified in $L/K$,
$
\operatorname{Frob}_{\mathfrak p,L/K}\subseteq C,
$
and
$$
N_{K/\mathbb Q}\mathfrak p
\le
\left(
4\log d_{\widetilde L}
+
2.5[\widetilde L:\mathbb Q]
+
5
\right)^2.
$$

\end{lem}

\begin{proof}
Applying Lemma \ref{prime ramify extn}, we see that every prime ideal in $S$ ramifies in the extension $M/K$. Since $L/K$ is Galois and $M/K$ is quadratic (hence Galois), the composite extension
$
\widetilde L=LM
$
is Galois over $K$.

First suppose that $M\subseteq L$. Then
$
\widetilde L=L.
$
Since every prime ideal in $S$ ramifies in $M/K$, and $M\subseteq L$, every prime ideal in $S$ ramifies in $L/K$.

By the effective Chebotarev Theorem (refer \ref{bachthm}), applied to the Galois extension $L/K$ and the conjugacy stable subset $C\subseteq \mathrm{Gal}(L/K)$, there exists a prime ideal $\mathfrak p\subseteq \mathcal O_K$ such that $\mathfrak p$ is unramified in $L/K$,
$
\operatorname{Frob}_{\mathfrak p,L/K}\subseteq C,
$
and
$$
N_{K/\mathbb Q}\mathfrak p
\le
\left(
4\log d_L+2.5[L:\mathbb Q]+5
\right)^2.
$$
Since $\widetilde L=L$, this gives
$$
N_{K/\mathbb Q}\mathfrak p
\le
\left(
4\log d_{\widetilde L}
+
2.5[\widetilde L:\mathbb Q]
+
5
\right)^2.
$$
Also, $\mathfrak p\notin S$, because every prime ideal in $S$ ramifies in $L/K$, whereas $\mathfrak p$ is unramified in $L/K$.

Now suppose that $M\not\subseteq L$. Since $M/K$ is quadratic, we have
$
L\cap M=K.
$
Hence  \cite[Theorem VI 1.14]{lang}
$$
\mathrm{Gal}(\widetilde L/K)
\simeq
\mathrm{Gal}(L/K)\times \mathrm{Gal}(M/K).
$$
Let $\operatorname{res}:\mathrm{Gal}(\widetilde L/K)\longrightarrow \mathrm{Gal}(L/K)
$
be the restriction map, and define
$
\widetilde C:=\operatorname{res}^{-1}(C).
$
Since $C$ is closed under conjugation in $\mathrm{Gal}(L/K)$, the set
$\widetilde C$ is closed under conjugation in $\mathrm{Gal}(\widetilde L/K)$. Applying the effective Chebotarev theorem to the Galois extension $\widetilde L/K$ and the conjugacy stable subset $\widetilde C$, we obtain a prime ideal $\mathfrak p\subseteq \mathcal O_K$ such that $\mathfrak p$ is unramified in $\widetilde L/K$,
$
\operatorname{Frob}_{\mathfrak p,\widetilde L/K}\subseteq \widetilde C,
$
and
$$
N_{K/\mathbb Q}\mathfrak p
\le
\left(
4\log d_{\widetilde L}
+
2.5[\widetilde L:\mathbb Q]
+
5
\right)^2.
$$

Since $\mathfrak p$ is unramified in $\widetilde L/K$, it is also unramified in $L/K$.
Moreover, every prime ideal in $S$ ramifies in $M/K$, hence also in $\widetilde L/K$. Therefore $\mathfrak p\notin S$.

Finally,
$
\operatorname{res}\left(\operatorname{Frob}_{\mathfrak p,\widetilde L/K}\right)
=
\operatorname{Frob}_{\mathfrak p,L/K}.
$
Since
$
\operatorname{Frob}_{\mathfrak p,\widetilde L/K}\subseteq \widetilde C
=
\operatorname{res}^{-1}(C),
$
we get
$
\operatorname{Frob}_{\mathfrak p,L/K}\subseteq C.
$
Thus, the required prime ideal $\mathfrak p$ exists.
\end{proof}

\subsection{Bound for the logarithm of the Discriminant}
In Lemma \ref{bach norm prime}, the norm of the prime ideal $\mathfrak p$ is bounded by $[\tilde{L}:\Q]$ and $\log d_{\tilde{L}}$. In our case, the first of these constants is absolutely bounded, and thus we need to bound the second. To do so, we improve the following theorem of Serre \cite[Proposition 6]{serre2} (a proof of which can be found in Mayle and Wang \cite[Lemma 7]{maylewang}) to Galois extensions of number fields. More specifically, the base field is now allowed to be any number field instead of $\mathbb Q$. Moreover, the lower bound has been tightened when $K=\Q$.
   \begin{lem}(Bound for Discriminant of Number Field)\label{bound disc} 
   Let $K$ be a number field and let $L/K$ be a nontrivial finite Galois extension. Put
$n=[L:K]$ and $f=[K:\mathbb Q]$. Let $d_L$ and $d_K$ denote the absolute discriminants of $L$ and $K$, respectively. Let $\mathfrak D_{L/K}\subseteq \mathcal O_L$ be the relative different ideal and let
$\mathfrak d_{L/K}= N_{L/K}(\mathfrak D_{L/K})\subseteq \mathcal O_K$
be the relative discriminant ideal. Define
$
\operatorname{rad}_K(\mathfrak d_{L/K})
=
\prod_{\mathfrak p\mid \mathfrak d_{L/K}}
N_{K/\mathbb Q}(\mathfrak p)$, then 
$$
n\log d_K+\frac n2\log \operatorname{rad}_K(\mathfrak d_{L/K})
\le
\log d_{L}\le
n\log d_K
+
(n-1)\log \operatorname{rad}_K(\mathfrak d_{L/K})
+
nf\log n.$$
 \end{lem}
\begin{proof}
Fix a nonzero prime ideal $\mathfrak{p}$ of $\mathcal{O}_{K}.$ In $\mathcal{O}_{L},$ let $\mathfrak{p}\mathcal{O}_{L}=\mathfrak{P}_{1}^{e_{1}}\cdots \mathfrak{P}_{g}^{e_{g}},$ where $\mathfrak{P}_{1}, \dots, \mathfrak{P}_{g}$ are distinct prime ideals of $\mathcal{O}_{L}$ lying above $\mathfrak{p}.$ For $\mathfrak{P}=\mathfrak{P}_{i},$  let $e_{\mathfrak P/\mathfrak p}$ denote its ramification index, and 
$f_{\mathfrak P/\mathfrak p}$ be its residue degree. Since $L/K$ is Galois, for fixed $\mathfrak p$, all primes $\mathfrak P\mid \mathfrak p$ have the same ramification index. Thus we may write
$e_{\mathfrak P/\mathfrak p}=e,$ and hence
$
\sum_{\mathfrak P\mid \mathfrak p}
e_{\mathfrak P/\mathfrak p}f_{\mathfrak P/\mathfrak p}
=
n.$
The relative discriminant formula \cite[Corollary 2.10, p. 202]{neukirch} gives
\begin{equation*}
  d_L=d_K^nN_{K/\mathbb Q}(\mathfrak d_{L/K}).   
\end{equation*}
Taking logarithms on both sides, we obtain
\begin{equation}\label{logdisc}
    \log d_L
=
n\log d_K+\log N_{K/\mathbb Q}(\mathfrak d_{L/K}).    
\end{equation}
 Using the factorization of the relative discriminant ideal as
$
\mathfrak d_{L/K}
=
\prod_{\mathfrak p}
\mathfrak p^{v_{\mathfrak p}(\mathfrak d_{L/K})}
$, we obtain
$
N_{K/\mathbb Q}(\mathfrak d_{L/K})
=
\prod_{\mathfrak p}
N_{K/\mathbb Q}(\mathfrak p)^{v_{\mathfrak p}(\mathfrak d_{L/K})}.
$
Again, by taking logarithms on both sides,
$$
\log N_{K/\mathbb Q}(\mathfrak d_{L/K})
=
\sum_{\mathfrak p}
v_{\mathfrak p}(\mathfrak d_{L/K})
\log N_{K/\mathbb Q}(\mathfrak p).
$$
Since
$
\mathfrak d_{L/K}=N_{L/K}(\mathfrak D_{L/K}),
$ we have $ v_{\mathfrak p}(\mathfrak d_{L/K})
=
\sum_{\mathfrak P\mid \mathfrak p}
f_{\mathfrak P/\mathfrak p}
v_{\mathfrak P}(\mathfrak D_{L/K}).
$ By Neukirch \cite[Theorem 2.6, p. 199]{neukirch}, we have $$ e_{\mathfrak{P}/\mathfrak{p}}-1\le v_{\mathfrak{P}}(\mathfrak{D}_{L/K})\le e_{\mathfrak{P}/\mathfrak{p}}-1+v_{\mathfrak{P}}(e_{\mathfrak{P}/\mathfrak{p}}). $$ Define $s_{\mathfrak{P}}:=v_{\mathfrak{P}}(\mathfrak{D}_{L/K})-(e_{\mathfrak{P}/\mathfrak{p}}-1).$ We thus have 
\begin{equation}\label{3.1valdisc}
 v_{\mathfrak P}(\mathfrak D_{L/K})
=
e_{\mathfrak P/\mathfrak p}-1+s_{\mathfrak P},
\quad \text{where}\quad 0\le s_{\mathfrak P}\le v_{\mathfrak P}(e_{\mathfrak P/\mathfrak p}).    
\end{equation}
 Fix a prime ideal $\mathfrak{p} \mid \mathfrak{d}_{L/ K}.$ By  \cite[Remark 3.39(b), p. 62]{milne2020algebraic}, we may assume that $\mathfrak p$ is ramified in $L/K$. Since $L/K$ is Galois, all primes above $\mathfrak p$ have the same ramification index $e$, and ramification means $e\ge 2$. Since $s_{\mathfrak P}\ge 0$, we get
$
v_{\mathfrak P}(\mathfrak D_{L/K})
\ge
e_{\mathfrak P/\mathfrak p}-1.
$
Therefore
$$
v_{\mathfrak p}(\mathfrak d_{L/K})
=
\sum_{\mathfrak P\mid \mathfrak p}
f_{\mathfrak P/\mathfrak p}
v_{\mathfrak P}(\mathfrak D_{L/K})
\ge
\sum_{\mathfrak P\mid \mathfrak p}
f_{\mathfrak P/\mathfrak p}
(e_{\mathfrak P/\mathfrak p}-1).
$$
Since $e_{\mathfrak P/\mathfrak p}=e$, this becomes
\begin{equation*}
 v_{\mathfrak p}(\mathfrak d_{L/K})
\ge
\sum_{\mathfrak P\mid \mathfrak p}
f_{\mathfrak P/\mathfrak p}(e-1).
\end{equation*}
The product of the local ramification indices and inertia degrees sum to the degree of the extension, that is, 
$
\sum_{\mathfrak P\mid \mathfrak p}
ef_{\mathfrak P/\mathfrak p}
=
n.
$
Therefore,
$
\sum_{\mathfrak P\mid \mathfrak p}
f_{\mathfrak P/\mathfrak p}
=
\frac ne,
$
hence
$
v_{\mathfrak p}(\mathfrak d_{L/K})
\ge
(e-1)\frac ne
=
n\left(1-\frac1e\right).$ Therefore
$v_{\mathfrak p}(\mathfrak d_{L/K})\ge \frac n2.$
Summing over all prime ideals $\mathfrak p\mid \mathfrak d_{L/K}$, we obtain
\begin{equation}\label{3.1}
    \log N_{K/\mathbb Q}(\mathfrak d_{L/K})
=
\sum_{\mathfrak p\mid \mathfrak d_{L/K}}
v_{\mathfrak p}(\mathfrak d_{L/K})
\log N_{K/\mathbb Q}(\mathfrak p)\ge
\sum_{\mathfrak p\mid \mathfrak d_{L/K}}
\frac n2
    \log N_{K/\mathbb Q}(\mathfrak p)=\dfrac{n}{2} \log \operatorname{rad}_{K}(\mathfrak{d}_{L/K}).
\end{equation}
Using the logarithm of the discriminant equation \eqref{logdisc} and combining with equation \eqref{3.1} gives us a lower bound
\begin{equation}\label{lowerbound}
   \boxed{ \log d_{L}\ge n\log d_{K}+\dfrac{n}{2}\log \operatorname{rad}_{K}(\mathfrak{d}_{L/K}).}
\end{equation}
Again use
$
v_{\mathfrak P}(\mathfrak D_{L/K})
=
e_{\mathfrak P/\mathfrak p}-1+s_{\mathfrak P}
$
with
$
0\le s_{\mathfrak P}\le v_{\mathfrak P}(e_{\mathfrak P/\mathfrak p}).
$ Then
\begin{equation*}
    v_{\mathfrak p}(\mathfrak d_{L/K})
=
\sum_{\mathfrak P\mid \mathfrak p}
f_{\mathfrak P/\mathfrak p}
v_{\mathfrak P}(\mathfrak D_{L/K}) \le
\sum_{\mathfrak P\mid \mathfrak p}
f_{\mathfrak P/\mathfrak p}
(e_{\mathfrak P/\mathfrak p}-1)
+
\sum_{\mathfrak P\mid \mathfrak p}
f_{\mathfrak P/\mathfrak p}
v_{\mathfrak P}(e_{\mathfrak P/\mathfrak p}).
\end{equation*}
Since $L/K$ is Galois, $e_{\mathfrak P/\mathfrak p}=e$ for all $\mathfrak P\mid \mathfrak p$. Hence
$$
\sum_{\mathfrak P\mid \mathfrak p}
f_{\mathfrak P/\mathfrak p}(e-1)
=
n-\sum_{\mathfrak P\mid \mathfrak p}
f_{\mathfrak P/\mathfrak p}\le n-1.
$$
For the second sum, by \cite[Equation 13, p. 66]{milne2020algebraic} we have $v_{\mathfrak{P}}(x)=e_{\mathfrak{P}/ \mathfrak{p}}v_{\mathfrak{p}}(x)$ for every $x\in K^{\times}.$ In particular, take $x=e$ then we have 
\begin{equation*}
v_{\mathfrak P}(e)
=
e_{\mathfrak P/\mathfrak p}v_{\mathfrak p}(e\mathcal O_K)
=
ev_{\mathfrak p}(e\mathcal O_K).
\end{equation*} 
Since $L/K$ is Galois, $e\mid n$. Hence
$
v_{\mathfrak p}(e\mathcal O_K)
\le
v_{\mathfrak p}(n\mathcal O_K).
$
Therefore
$
v_{\mathfrak P}(e)
\le
ev_{\mathfrak p}(n\mathcal O_K).
$
Thus
$$
\sum_{\mathfrak P\mid \mathfrak p}
f_{\mathfrak P/\mathfrak p}
v_{\mathfrak P}(e)
\le
\sum_{\mathfrak P\mid \mathfrak p}
f_{\mathfrak P/\mathfrak p}
ev_{\mathfrak p}(n\mathcal O_K).
$$
Using
$
\sum_{\mathfrak P\mid \mathfrak p}
ef_{\mathfrak P/\mathfrak p}=n,
$
we obtain
$
\sum_{\mathfrak P\mid \mathfrak p}
f_{\mathfrak P/\mathfrak p}
v_{\mathfrak P}(e)
\le
nv_{\mathfrak p}(n\mathcal O_K).
$
Consequently,
$$
v_{\mathfrak p}(\mathfrak d_{L/K})
\le
n-1+nv_{\mathfrak p}(n\mathcal O_K).
$$
Multiplying by $\log N_{K/\mathbb Q}(\mathfrak p)$ and summing over all $\mathfrak p\mid \mathfrak d_{L/K}$ gives
$$
\log N_{K/\mathbb Q}(\mathfrak d_{L/K})
\le
\sum_{\mathfrak p\mid \mathfrak d_{L/K}}
(n-1)\log N_{K/\mathbb Q}(\mathfrak p)
+
n\sum_{\mathfrak p\mid \mathfrak d_{L/K}}
v_{\mathfrak p}(n\mathcal O_K)
\log N_{K/\mathbb Q}(\mathfrak p).
$$
The first sum is
$$
(n-1)\log \operatorname{rad}_K(\mathfrak d_{L/K}).
$$
For the second sum, we enlarge the sum to all prime ideals of $\mathcal O_K$ (the sum is finite as $v_{\mathfrak p}(n\mathcal O_K)=0$ for almost all $\mathfrak p$):
$$
\sum_{\mathfrak p\mid \mathfrak d_{L/K}}
v_{\mathfrak p}(n\mathcal O_K)
\log N_{K/\mathbb Q}(\mathfrak p)
\le
\sum_{\mathfrak p}
v_{\mathfrak p}(n\mathcal O_K)
\log N_{K/\mathbb Q}(\mathfrak p).
$$
Writing
$
n\mathcal O_K
=
\prod_{\mathfrak p}
\mathfrak p^{v_{\mathfrak p}(n\mathcal O_K)},
$
and taking logarithms on both sides, we have
\begin{equation*}
 \sum_{\mathfrak p}
v_{\mathfrak p}(n\mathcal O_K)
\log N_{K/\mathbb Q}(\mathfrak p)
=
\log N_{K/\mathbb Q}(n\mathcal O_K).   
\end{equation*}
But
$
N_{K/\mathbb Q}(n\mathcal O_K)=n^f.
$
Hence
$
\log N_{K/\mathbb Q}(n\mathcal O_K)=f\log n.
$
Therefore
$$
\log N_{K/\mathbb Q}(\mathfrak d_{L/K})
\le
(n-1)\log \operatorname{rad}_K(\mathfrak d_{L/K})
+
nf\log n.
$$
Using this with equation \eqref{logdisc} gives us an upper bound
\begin{equation}\label{upperbound}
   \boxed{ \log d_{L}\le
n\log d_K
+
(n-1)\log \operatorname{rad}_K(\mathfrak d_{L/K})
+
nf\log n.}
\end{equation}
Combining the lower bound \eqref{lowerbound} and upper bound \eqref{upperbound} we obtain the desired inequality.
\end{proof}

\section{An Effective version of Faltings Isogeny Theorem for Number Fields}
In this section, we prove an Effective version of Faltings isogeny theorem for number fields. This theorem lies at the heart of providing a bound on the primes for which the image cannot be contained in the normalizer of a non-split Cartan subgroup. We begin by describing the primes of bad reduction for an abelian variety $A$ over a number field $K$ in terms of the primes that ramify in some extension in $\S$\ref{S:primes of bad reduction}. We then prove some number theoretic lemmas in $\S$\ref{S:some number theory} which will be used to prove the effective isogeny theorem in $\S$\ref{S: effective isogeny theorem for number fields}.
\subsection{Primes of bad reduction}\label{S:primes of bad reduction}
Let $A$ be an abelian variety over the number field $K$ and let $K(A[m])$ be the field obtained by adjoining the co-ordinates of the $m$-torsion points of $A$ in $A(\overline{\Q})$ to the field $K$. For a prime number $\ell$, we write
\begin{equation*}
    K(A[\ell^\infty])=\bigcup_{n=1}^{\infty} K(A[\ell^n]).
\end{equation*}
The following theorem by Serre and Tate relates the primes of good reduction of $A$ with the primes that ramify in the extension $K(A[\ell^\infty])$.
\begin{theorem}[Serre-Tate]\label{serretate}
    Let $A/K$ be an abelian variety. Let $\mathfrak p$ denote a prime ideal in $\mathcal{O}_K$ and let $p=\text{char}(\mathcal{O}_K/\mathfrak p\mathcal{O}_K)$. The following are equivalent:
    \begin{enumerate}
        \item $A$ has good reduction at $\mathfrak p$.
        \item $K(A[m])/K$ is unramified at $\mathfrak p$ for all $m$ prime to $p$.
        \item There exist infinitely many integers $m$, prime to $p$, such that $K(A[m]/K$ is unramified at $\mathfrak p$.
    \end{enumerate}
\end{theorem}
\begin{proof}
    \cite[Theorem 1]{serretate}
\end{proof}
\begin{corollary}\label{ram primes}
    Let $B_A$ be the product of ideal $(\ell)$ and the prime ideals of bad reduction for $A$. The extension $K(A[\ell^\infty])/K$ is ramified at exactly the prime ideals dividing $B_A$.
\end{corollary}
\begin{proof}
    Suppose that $\mathfrak q$ is a prime ideal lying above the rational prime $\ell$, that is, $\mathfrak q\mid (\ell)$ as ideals. We show that $\mathfrak q$ ramifies in $K(A[\ell^\infty])/K$. For a finite extension $L$ of $K$, the ramification of the prime ideal $\mathfrak q$ with respect to a prime ideal $\mathfrak Q$ of $L$, is controlled by the ramification of the extension of the completions $L_\mathfrak{Q}/K_\mathfrak{q}$. To be more precise,
    \begin{equation}\label{ram local global}
        e(\mathfrak Q\mid\mathfrak q)= e(L_{\mathfrak Q}, \mid K_\mathfrak{q}),
    \end{equation}
    where the left-hand side denotes the ramification index of the prime ideal $\mathfrak q$ with respect to $\mathfrak Q$ and the right-hand side denotes the ramification index of the extension of completions.\\
    Let $\zeta_n$ denote a primitive $n$-th root of unity and let $L_n=K(\zeta_{\ell^n})$. For $n$ large enough, $L_n$ is a non-trivial finite extension of $K$. Let $K_\mathfrak{q}$ denote the completion of $K$ with respect to the prime ideal $\mathfrak q$. This is a finite extension of $\mathbb Q_\ell$. Let $E_n=\mathbb Q_\ell(\zeta_{\ell^n})$ and let $M_n=K_\mathfrak q(\zeta_{\ell^n})$. Pictorially, we have    
    \begin{equation*}
    \begin{tikzcd}
                 & M_n \arrow[ld, dash] \arrow[rd, dash] &              \\
    K_\mathfrak q \arrow[rd, dash] &                         & E_n \arrow[ld, dash] \\
                 & \mathbb{Q}_\ell                       &             
    \end{tikzcd}
    \end{equation*}
    Assume that $M_n$ over $K_\mathfrak{q}$ is an unramified extension for all integers $n$. Then 
    \begin{equation}\label{ram indices}
        e(K_\mathfrak{q}\mid \Q_\ell)= e(M_n\mid \Q_\ell).
    \end{equation}
    By the tower law for ramification indices, $e(E_n\mid \Q_\ell)$ divides $e(M_n\mid \Q_\ell)$. The extension $E_n$ over $\Q_\ell$ is totally ramified of degree $\varphi(\ell^n)$, where $\varphi$ denotes the Euler totient function. Using equation (\ref{ram indices}), we get
    \begin{equation*}
        \varphi(\ell^n)\mid e(K_\mathfrak{q}\mid \Q_\ell)\leq [K_\mathfrak{q}:\Q_\ell].
    \end{equation*}
    This leads to a contradiction, as for large enough $n$ this inequality does not hold, since the right-hand side is finite. Thus $M_n$ over $K_\mathfrak{q}$ is ramified for some $n$ and by equation (\ref{ram local global}), we see that $\mathfrak{q}$ is ramified in $L_n$. By the Weil Pairing (see exercise in \cite[p. 55]{serremordellweil}) on the abelian variety A, we have the following tower of fields
    \begin{equation*}
        K\subseteq K(\zeta_{\ell^n})\subseteq K(A[\ell^n])\subseteq K(A[\ell^\infty]).
    \end{equation*}
    Hence, $\mathfrak{q}$ is ramified in $K(A[\ell^\infty])$.\\
    The proof for prime ideals $\mathfrak p$ not dividing the ideal $(\ell)$ is the same as in Mayle and Wang \cite[Corollary 9, p. 1406]{maylewang} with a few minor changes. That is, we replace the rational prime $p$ by the prime ideal $\mathfrak p$ and $\Q$ by the number field $K$. Moreover, we use the general version of the Serre-Tate Theorem \ref{serretate}. 
\end{proof}

\begin{remark}\label{R: bad primes}
    Let $A=E_1\times E_2$ be the abelian variety that is the product of the two elliptic curves over $K$ with conductors $\mathfrak N_{E_1}$ and $\mathfrak N_{E_2}$ respectively. Theorem \ref{serretate} states that $A$ has good reduction at a prime $\mathfrak p$ lying over $\ell$ if and only if both $E_1$ and $E_2$ have good reduction at $\mathfrak p$. Thus $\textrm{rad}(B_A)=\textrm{rad}(\ell\mathcal{O}_K\mathfrak{N}_{E_1}\mathfrak{N}_{E_2})$.
\end{remark}
\subsection{Trace of the Frobenius}
Let $E$ be an elliptic curve over a number field $K$. For each rational prime $p$, and any prime ideal $\mathfrak p$ lying over $p$, we may associate a number 
\begin{equation}\label{trace frob}
    a_{\mathfrak p}(E):=N(\mathfrak p)+1-\#E(\mathcal{O}_K/\mathfrak p),
\end{equation}
where $N(E)$ denotes the norm of the prime ideal $\mathfrak p$, and $\#E(\mathcal{O}_K/\mathfrak p)$ denote the number of solutions for $E$ in $\mathcal{O}_K/\mathfrak p$. This definition seems to depend on the choice of prime ideal $\mathfrak p$ lying over $p$. We will now show that this is not the case when $K$ is Galois over $\Q$.
  \begin{propn}
     Let $K$ be a Galois extension of $\Q$ and let $E$ be an elliptic curve over $K$. Let $\mathfrak p$ and $\mathfrak q$ be two prime ideals lying over the rational prime $p$. We have $a_\mathfrak p(E)=a_\mathfrak q(E)$. 
  \end{propn}
  \begin{proof}
      In a Galois extension, all prime ideals lying above a rational prime are conjugate to each other. That is, we can find a $\sigma\in\mathrm{Gal}(K/\Q)$ such that $\sigma(\mathfrak p)=\mathfrak q$. It follows that the norms are the same and the residue fields are isomorphic(see \cite[Chapter 1, \S 9]{neukirch}) Looking at the terms appearing in equation (\ref{trace frob}), we get the desired result.
  \end{proof}

Faltings isogeny theorem for number fields states that isogenous elliptic curves $E$ and $E'$  have $a_\mathfrak p(E)=a_\mathfrak p(E')$ for all primes of good reduction $\mathfrak p$. Hence, if two elliptic curves are not isogenous, there exists a prime ideal $\mathfrak p$ for which the traces are unequal. We give a bound on the norm of one such prime ideal. This bound will then be used to get a value of $c(E)$. The main goal of this section is to prove an effective version of Faltings isogeny theorem for number fields. To do so, we first establish several number-theoretic lemmas.\\
 \subsection{Some number theory}\label{S:some number theory}
 
 We first produce an integral basis using the embedding of the number field $K$ in Euclidean space by means of the  Minkowski map. 
 \begin{lem}\label{Important lem}
     Let $K$ be a number field of degree $[K:\mathbb{Q}]=n_{K},$ and $d_{K}=\left|\operatorname{Disc}(K/\mathbb{Q})\right|$ denote the absolute discriminant of $K.$ Then $\mathcal{O}_{K}$ admits an integral basis $w_{1},\dots, w_{n_{K}}$ such that $|\sigma(w_{j})|\le n_{K} \sqrt{d_{K}}$ for every embedding $\sigma: K\hookrightarrow \mathbb{C}$ and for every $j=1,2,\dots,n_{K}.$ Moreover, $$ \prod\limits_{\sigma:K\hookrightarrow \mathbb{C}}\left(\sum\limits_{j=1}^{n_{K}}|\sigma (w_{j})|\right)\le n_{K}^{2n_{K}}d_{K}^{n_{K}/2}.$$
     \end{lem}
 \begin{proof}
     For $n_{K}=1,$ take the basis $\{1\}.$ Now consider $n_{K}\ge 2.$ Let $r_{1}$ be the number of real embeddings of $K,$ and let $r_{2}$ be the number of conjugate pairs of nonreal embeddings. Thus $n_{K}=r_{1}+2r_{2}$. Choose one embedding from each complex conjugate pair and consider the embedding $$ \iota:K\longrightarrow \mathbb{R}^{r_{1}}\times \mathbb{C}^{r_{2}},\ \alpha\longmapsto (\sigma_{1}(\alpha),\dots,\sigma_{r_{1}}(\alpha),\tau_{1}(\alpha),\dots,\tau_{r_{2}}(\alpha)).$$ Now apply \cite[Proposition 4.26, p. 79]{milne2020algebraic} to the unit ideal $\mathcal{O}_{K}.$ Since $N_{K/\mathbb{Q}}(\mathcal{O}_{K})=1,$ the image $\iota(\mathcal{O}_{K}):=\Lambda$ is a full lattice of covolume $2^{-r_{2}}\sqrt{d_{K}}.$ Consider the set $B=[-1,1]^{r_{1}}\times D^{r_{2}},$ where $D$ is a closed unit disc in $ \mathbb{C}. $ Then it is compact, convex, symmetric about the origin, and has nonempty interior. Under the identification $\mathbb{R}^{r_{1}}\times \mathbb{C}^{r_{2}}\cong \mathbb{R}^{n_{K}}$ via $$ (x_{1},\dots,x_{r_{1}},z_{1},\dots,z_{r_{2}})\longmapsto (x_{1},\dots,x_{r_{1}}, \operatorname{Re}z_{1}, \operatorname{Im}z_{1},\dots, \operatorname{Re}z_{r_{2}},\operatorname{Im}z_{r_{2}}),$$ let $\operatorname{vol}$ denote $n_{K}$-dimensional Lebesgue measure. By repeatedly applying \cite[Theorem 2.37(a), p. 67]{folland1999real} the volume of $B$ is the product of the lengths of its interval factors and the areas of its disc factors, hence $\operatorname{vol}(B)=2^{r_{1}}\pi^{r_{2}}.$ Let $\lambda_{1}\le \dots \le \lambda_{n_{K}}$ be its successive minima with respect to $\Lambda.$ Here $\lambda_{i}$ is the smallest scaling factor such that $\lambda_{i}B$  contains at least $i$ linearly independent lattice vectors. First, we show that $\lambda_{i}\ge 1 \ (i=1,\dots,n_{K}).$ Indeed, if there exists nonzero $\alpha\in \mathcal{O}_{K}$ such that $\iota(\alpha)\in tB$ for some $0<t<1,$ where $tB=[-t,t]^{r_{1}}\times (tD)^{r_{2}}.$ This implies that $\operatorname{vol}(tB)=2^{r_{1}}\pi^{r_{2}}t^{n_{K}}.$ Moreover, $|\sigma_{i}(\alpha)|\le t$ for all $i=1,2\dots,r_{1}$ and $|\tau_{j}(\alpha)|\le t$ for all $j=1,\dots,r_{2}.$ Then $$ |N_{K/\mathbb{Q}}(\alpha)|=\prod\limits_{i=1}^{r_{1}}|\sigma_{i}(\alpha)|\prod\limits_{j=1}^{r_{2}}|\tau_{j}(\alpha)|^{2}\le t^{r_{1}}t^{2r_{2}}=t^{n_{K}}<1,$$ which is a contradiction. Hence, by the definition of the first successive minimum, $\lambda_{1}\ge 1$ and thus $\lambda_{i}\ge 1$ for all $i$. Now apply Minkowski's second body theorem \cite[Theorem 2.9, 2nd Chapter, p. 24]{Evertse} we obtain $$ \lambda_{1}\cdots\lambda_{n_{K}}\le \dfrac{2^{n_{K}}\operatorname{det}\Lambda}{\operatorname{vol}(B)}=\dfrac{2^{n_{K}}\cdot 2^{-r_{2}}\sqrt{d_{K}}}{2^{r_{1}}\pi^{r_{2}}}=\left(\dfrac{2}{\pi}\right)^{r_{2}}\sqrt{d_{K}}\le \sqrt{d_{K}}.$$ Since all $\lambda_{1},\dots, \lambda_{n_{K}-1}\ge 1,$ so we have $\lambda_{n_{K}}\le \sqrt{d_{K}}.$ Choose linearly independent lattice vectors $v_{1}, \dots, v_{n_{K}}$ with $v_{i}\in \lambda_{i}B.$ Then by applying \cite[Lemma 2.1, p. 12]{Evertse} we get a basis $b_{1},\dots,b_{n_{K}}$ of $\Lambda$ of the form $b_{i}=\sum\limits_{j=1}^{n_{K}}a_{ij}v_{j},\ 0\le a_{ij}\le 1.$ The norm whose unit ball is $B$ is the maximum of the absolute values of the real and complex coordinates. By the triangle inequality, $$ \left\|b_{i}\right\|_{B}\le \sum\limits_{j=1}^{n_{K}}a_{ij}\left\|v_{j}\right\|_{B}\le \sum\limits_{j=1}^{n_{K}}\lambda_{j}\le n_{K}\lambda_{n_{K}}\le n_{K}\sqrt{d_{K}}.$$ Write $b_{i}=\iota(w_{i})$ for some $w_{i}\in \mathcal{O}_{K}.$ Since $\iota$ is injective and $b_{i}$'s form a lattice basis of $\iota(\mathcal{O}_{K}),$ the elements $w_{1},\dots,w_{n_{K}}$ form an integral basis of $\mathcal{O}_{K}.$ The above estimation says $|\sigma(w_{i})|\le n_{K}\sqrt{d_{K}}$ for all embeddings $\sigma: K\hookrightarrow \mathbb{C}$ and for all $i=1,\dots,n_{K}.$ Consequently, for each embedding $\sigma,$ $\sum\limits_{j=1}^{n_{K}}|\sigma(w_{j})|\le n_{K}^{2}\sqrt{d_{K}}.$ Thus $\prod\limits_{\sigma:K\hookrightarrow \mathbb{C}}\left(\sum\limits_{j=1}^{n_{K}}|\sigma (w_{j})|\right)\le n_{K}^{2n_{K}}d_{K}^{n_{K}/2}.$
 \end{proof}
 Using the special basis obtained in Lemma \ref{Important lem}, we produce an element of the number field whose norm satisfies a particular bound. This bound depends only on the number field $K$ and the norms of a prescribed set of prime ideals. 

 \begin{lem}\label{norm bound}
     Let $K$ be a number field of degree $[K:Q]=n_{K},$ $d_{K}=|\operatorname{Disc(K/\mathbb{Q)|}}$ and let $\mathfrak{m}=\mathfrak{q}_{1}\cdots \mathfrak{q}_{e}, \ e\ge 1,$ where $\mathfrak{q}_{1},\dots, \mathfrak{q}_{e}$ are distinct nonzero prime ideals of $\mathcal{O}_{K}.$ Then there exists $\alpha\in \mathcal{O}_{K}\setminus\{0\}$ such that $v_{\mathfrak{q}_{i}}(\alpha)=1\ \text{for}\ 1\leq i\leq e$ and $|N_{K/\mathbb{Q}}(\alpha)|\le n_{K}^{2n_{K}}d_{K}^{n_{K}/2}(N_{K/\mathbb{Q}}\mathfrak{m})^{2n_{K}}.$ 
     \end{lem}
     \begin{proof}
         By Lemma \ref{Important lem}, choose an integral basis $w_{1},\dots,w_{n_{K}}$ of $\mathcal{O}_{K}$ such that $$\prod\limits_{\sigma:K\hookrightarrow \mathbb{C}}\left(\sum\limits_{j=1}^{n_{K}}|\sigma (w_{j})|\right)\le n_{K}^{2n_{K}}d_{K}^{n_{K}/2}.$$ Since $\mathcal{O}_{K}$ is a Dedekind domain, the unique factorization of nonzero ideals implies that $\mathfrak{q}_{i}^{2}\subsetneq \mathfrak{q}_{i}.$ Choose $a_{i}\in \mathfrak{q}_{i}\setminus \mathfrak{q}_{i}^{2},$ then by Chinese remainder theorem, there exists $\alpha_{0}\in \mathcal{O}_{K}$ such that $\alpha_{0}\equiv a_{i} \pmod{\mathfrak{q}_{i}^{2}}\quad (i=1,\dots,e).$ Put $u=N_{K/\mathbb{Q}}(\mathfrak{m}^{2})=(N_{K/\mathbb{Q}}\mathfrak{m})^{2}.$ Since the additive group $\mathcal{O}_{K}/\mathfrak{m}^{2}$ has order $u,$ we have $u\mathcal{O}_{K}\subseteq \mathfrak{m}^{2}.$ Now write $\alpha_{0}=\sum\limits _{j=1}^{n_{K}}b_{j}w_{j},\ b_{j}\in \mathbb{Z},$ then for each $j,$ the division algorithm gives integers $t_{j},c_{j}$ such that $b_{j}=ut_{j}+c_{j},\ 0\le c_{j}<u.$ Define $\alpha=\sum\limits_{j=1}^{n_{K}}c_{j}w_{j},$ then $$\alpha-\alpha_{0}=\sum\limits_{j=1}^{n_{K}}(c_{j}-b_{j})w_{j}=-u\sum\limits_{j=1}^{n_{K}}t_{j}w_{j}.$$ Therefore, $$ \alpha-\alpha_{0}\in u \mathcal{O}_{K}\subseteq\mathfrak{m}^{2}=\mathfrak{q}_{1}^{2}\cdots \mathfrak{q}_{e}^{2}\subseteq \mathfrak{q}_{i}^{2}.$$ Thus $\alpha-a_{i}=(\alpha-\alpha_{0})+(\alpha_{0}-a_{i})$ belongs to $\mathfrak{q}_{i}^{2}.$ Hence $\alpha \equiv a_{i}\pmod{\mathfrak{q}_{i}^{2}}\ (i=1,\dots,e).$ Thus $\alpha=a_{i}+(\alpha-a_{i})\in \mathfrak{q}_{i}.$  But $\alpha \notin \mathfrak{q}_{i}^{2}.$ Suppose not, then $a_{i}=\alpha-(\alpha-a_{i})\in \mathfrak{q}_{i}^{2},$ contrary to our choice. Therefore $\alpha \in \mathfrak{q}_{i}\setminus \mathfrak{q}_{i}^{2},$ so $\alpha\ne 0$ and $v_{\mathfrak{q}_{i}}(\alpha)=1.$ Let $\sigma: K\hookrightarrow \mathbb{C}$ be an embedding, then $\sigma(\alpha)=\sum\limits_{j=1}^{n_{K}}c_{j}\sigma(w_{j})$ and hence $$ |\sigma(\alpha)|\le \sum \limits_{j=1}^{n_{K}}c_{j}|\sigma (w_{j})|\le u \sum \limits_{j=1}^{n_{K}}|\sigma (w_{j})|.$$ Now \begin{align*}
             |N_{K/\mathbb{Q}}(\alpha)|&=\prod\limits_{\sigma:K\hookrightarrow \mathbb{C}}|\sigma(\alpha)|
             \le \prod\limits_{\sigma:K\hookrightarrow \mathbb{C}}\left(u \sum \limits_{j=1}^{n_{K}}|\sigma(w_{j})|\right)\\
             &\le u^{n_{K}} n_{K}^{2n_{K}}d_{K}^{n_{K}/2}= n_{K}^{2n_{K}}d_{K}^{n_{K}/2}(N_{K/\mathbb{Q}}\mathfrak{m})^{2n_{K}}.
         \end{align*}
     \end{proof}
 \begin{remark}
     Note that the $\alpha$ obtained in Lemma \ref{norm bound} is compatible with Lemma \ref{prime ramify extn} and Lemma \ref{bach norm prime}.
 \end{remark} 
 \subsection{Effective Isogeny Theorem over Galois Number Fields} \label{S: effective isogeny theorem for number fields}
  We use the element produced in Lemma \ref{norm bound}, in conjunction with Lemma \ref{prime ramify extn} to produce a field extension in which a prescribed set of prime ideals ramify. Finally, we use the modified Bach and Sorenson Theorem \ref{bach norm prime} to get a prime ideal outside the prescribed set. The prime $2$ appears in the following theorem as we used the $2$-adic representation associated to the elliptic curve to get the bound in the first line of the proof. We couple this with Remark \ref{R: bad primes} to exclude the primes of bad reduction.
  \begin{theorem}[Effective Isogeny Theorem over Number Fields]\label{Eff iso for K}
Assume GRH. Let $K$ be a Galois number field of degree  $n_{K}$ over $\mathbb{Q}$, and let $E_1,E_2$ be elliptic curves without complex multiplication defined over $K$. Suppose that $E_1$ and $E_2$ are not $K$-isogenous. Then there exists a prime ideal $\mathfrak p\subseteq \mathcal O_K$ of good reduction for both $E_1$ and $E_2$ such that
$$
a_{\mathfrak p}(E_1)\neq a_{\mathfrak p}(E_2),
$$
and satisfying
$$
N_{K/\mathbb Q}\mathfrak p
\le
\left(
C_1(K)
\log N_{K/\mathbb Q}
\operatorname{rad}
\left(
2\mathcal O_K\mathfrak N_{E_1}\mathfrak N_{E_2}
\right)
+
C_2(K)
\right)^2,
$$
where $\mathfrak N_{E_i}$ denotes the conductor ideal of $E_i$, and $C_1(K),C_2(K)$ are constants depending only on $K$.
The constants are explicitly given as
\begin{align*}
    C_1(K)=& 510.352+513.024n_{K},\\
    C_2(K)=& 513.024\bigl(1+\frac{n_{K}}{4}\log d_K
+n_{K}\log(384n_{K})\bigr)
+1324.8n_{K}+45.7,
\end{align*}
where $d_K$ denotes the absolute discriminant of the number field $K.$
\end{theorem}
\begin{proof}
We follow the same proof as in Mayle and Wang to conclude that $[L:K]\leq 192$ \cite[Proof of Theorem 2, p. 1408]{maylewang}.

Define
$$
\mathfrak m
:=
\operatorname{rad}
\left(
2\mathcal O_K\mathfrak N_{E_1}\mathfrak N_{E_2}
\right)
=
\prod_{\mathfrak q\mid 2\mathcal O_K\mathfrak N_{E_1}\mathfrak N_{E_2}}
\mathfrak q,
$$ and 
let
$
S:=\{\mathfrak q\subseteq \mathcal O_K:\mathfrak q\mid \mathfrak m\}.
$ By Lemma \ref{norm bound}, choose $\alpha\in \mathcal{O}_{K}\setminus \{0\}$ satisfying 
$
v_{\mathfrak q}(\alpha)= 1
$
for every $\mathfrak q\in S,$ and 
\begin{equation}\label{4.2 log norm}
    \log|N_{K/\mathbb{Q}}(\alpha)|\le 2n_{K}\log n_{K}+ \dfrac{n_{K}}{2} \log d_{K}+ 2n_{K}\log N_{K/\mathbb{Q}}\mathfrak{m}.
\end{equation}
Fix this $\alpha$ and put
$
M=K(\sqrt{\alpha})$ and $ \widetilde{L}=LM. 
$ Since $v_{\mathfrak{q}}=1$ for every $\mathfrak{q}\in S,$  Lemma \ref{prime ramify extn} implies that every prime ideal $\mathfrak q\in S$ ramifies in $M/K$ .
 Applying Lemma \ref{bach norm prime} with the same $\alpha$, we obtain  a prime ideal $\mathfrak p\subseteq \mathcal O_K$ such that
$
\mathfrak p\notin S,
$
$\mathfrak p$ is unramified in $L/K$,
$
\operatorname{Frob}_{\mathfrak p,L/K}\subseteq C,
$
and
$
N_{K/\mathbb Q}\mathfrak p
\leq
\left(
4\log d_{\widetilde L}
+
2.5[\widetilde L:\mathbb Q]
+
5
\right)^2.
$

Since $\mathfrak p\notin S$, we have
$
\mathfrak p\nmid 2\mathcal O_K\mathfrak N_{E_1}\mathfrak N_{E_2}.
$
Therefore, $\mathfrak p$ is not above $2$, and both $E_1$ and $E_2$ have good reduction at $\mathfrak p$. So again by \cite[Proposition 12, part (3)]{maylewang}, we have
$$
\operatorname{tr}\rho_1(\operatorname{Frob}_{\mathfrak p})
\neq
\operatorname{tr}\rho_2(\operatorname{Frob}_{\mathfrak p}).
$$

Hence
$
a_{\mathfrak p}(E_1)\neq a_{\mathfrak p}(E_2).
$ It remains to prove the required bound. Furthermore,
$
[\widetilde L:K]\leq 2[L:K]\leq 384.
$
Therefore,
$
[\widetilde L:\mathbb Q]
=
[\widetilde L:K]n_{K}
\leq
384n_{K}.
$ Thus, if
$
\log d_{\widetilde L}\leq 100,
$
then by applying the effective Chebotarev bound \ref{bachthm} with the constants $4,2.5,$ and $5$, we find that
$$
N_{K/\mathbb Q}\mathfrak p
\leq
\left(
4\cdot 100+2.5\cdot 384n_{K}+5
\right)^2
=
\left(
405+960n_{K}
\right)^2.
$$

If
$
100\leq \log d_{\widetilde L}\leq 1000,
$ taking the maximum of each coefficient over the degree columns in the corresponding row of \cite[Table 1, p. 1730]{bach1996explicit}, we obtain 
\begin{align*}
    N_{K/\mathbb{Q}}\mathfrak{p} &\le (1.796 \cdot 1000 + 1.89 \cdot 384 n_{K} + 6.3 )^{2}\\
    &=(1796+725.76n_{K} + 6.3)^{2}\\
    &=(1802.3 + 725.76 n_{K})^{2}.
\end{align*} If $\log d_{\widetilde{L}}\ge 1000,$ then again using \cite[Table 1, p. 1730]{bach1996explicit}, we obtain \begin{align*}
    N_{K/\mathbb{Q}}\mathfrak{p}&\le (1.336 \log d_{\widetilde{L}} + 3.45 \cdot 384 n_{K} + 45.7 )^{2}\\
    &=(1.336 \log d_{\widetilde{L}}+1324.8 n_{K}+ 45.7)^{2}.
\end{align*}
Therefore, in all cases, 
\begin{equation}\label{4.3 norm prime ideal}
    \begin{aligned}
    N_{K/\mathbb{Q}}\mathfrak{p}\le\max\bigl\{&\bigl(405+960n_{K}\bigr)^{2},\\ 
    & \bigl(1802.3+725.76n_{K}\bigr)^{2},\\
    & \bigl(1.336 \log d_{\widetilde{L}}+ 1324.8 n_{K}+45.7\bigr)^{2}\bigr\}.
\end{aligned}
\end{equation}
Now it remains to bound $\log d_{\widetilde{L}}.$ From Lemma \ref{bound disc}, we have $$ \log d_{L}\le [L:K] \log d_{K}+([L:K]-1)\log \operatorname{rad}_{K}(\mathfrak{d}_{L/K})+[L:K]n_{K}\log[L:K].$$ Recall that $\operatorname{rad}_{K}(\mathfrak{d}_{L/K})=\prod\limits_{\mathfrak{q}\mid \mathfrak{d}_{L/K}}N_{K/\mathbb{Q}}\mathfrak{q}. $ A prime ideal divides the relative discriminant precisely when it ramifies in $L/K$ \cite[Remark 3.39(b), p. 62]{milne2020algebraic}. Since $L/K$ is unramified outside $S,$ every prime ideal dividing $\mathfrak{d}_{L/K}$ belongs to $S.$ Therefore 
\begin{align*}
   \operatorname{rad}_{K}(\mathfrak{d}_{L/K})\le \prod\limits_{\mathfrak{q}\in S}N_{K/\mathbb{Q}} \mathfrak{q}=N_{K/\mathbb{Q}}(\mathfrak{q}_{1}\cdots \mathfrak{q}_{e} )=N_{K/\mathbb{Q}}\mathfrak{m}. 
\end{align*}
Hence $\log \operatorname{rad}_{K}(\mathfrak{d}_{L/K})\le \log N_{K/\mathbb{Q}}\mathfrak{m}.$ Consequently, 
\begin{equation}\label{4.4 log disc}
     \log d_{L}\le 192 \log d_{K}+191\log N_{K/\mathbb{Q}}\mathfrak{m}+192 n_{K}\log 192.
\end{equation}
Suppose first that $[\widetilde{L}:L]=2,$ i.e. $\sqrt{\alpha}\notin L.$ Then the discriminant of the $L$-basis $\{1,\sqrt{\alpha}\}$ is $4\alpha.$ So the relative field discriminant divides the discriminant of this order, i.e.,  $\mathfrak{d}_{\widetilde{L}/L}\mid (4\alpha)\mathcal{O}_{L},$ which implies $(4\alpha)\mathcal{O}_{L}=\mathfrak{d}_{\widetilde{L}/L}\mathfrak{a}$ for some nonzero integral ideal $\mathfrak{a}$ of $\mathcal{O}_{L}.$ Now apply \cite[Theorem 6.1, p. 11]{rabinoff4803discriminants} to $\mathbb{Q}\subseteq L \subseteq \widetilde{L}.$ Since $[\widetilde{L}:L]=2,$ it gives $d_{\widetilde{L}}=d_{L}^{2}N_{L/\mathbb{Q}}(\mathfrak{d}_{\widetilde{L}/L})$. Consequently, $N_{L/\mathbb{Q}}((4\alpha)\mathcal{O}_{L})=N_{L/\mathbb{Q}}(\mathfrak{d}_{\widetilde{L}/L})N_{L/\mathbb{Q}}(\mathfrak{a}).$ Thus $$ N_{L/\mathbb{Q}}(\mathfrak{d}_{\widetilde{L}/L})\le N_{L/\mathbb{Q}}((4\alpha)\mathcal{O}_{L})=|N_{L/\mathbb{Q}}(4\alpha)|.$$ The above follows from \cite[Proposition 4.1, Proposition 4.2(a), p. 68-89]{milne2020algebraic}. So that 
\begin{equation}\label{4.5 disc comp}
    d_{\widetilde{L}}\le d_{L}^{2}|N_{L/\mathbb{Q}}(4\alpha)|.
\end{equation}
Now by transitivity of field norms, for $\mathbb{Q}\subseteq K \subseteq L,$ we have  
\begin{align*}
    N_{L/\mathbb{Q}}(4\alpha)&=N_{K/\mathbb{Q}}(N_{L/ K}(4\alpha))=N_{K/\mathbb{Q}}((4\alpha)^{[L:K]})\\
    &=(N_{K/ \mathbb{Q}}(4\alpha))^{[L:K]}= (4^{n_{K}} N_{K/\mathbb{Q}}(\alpha))^{[L:K]}
\end{align*}
Now, from equation \eqref{4.5 disc comp}, we get 
\begin{equation}\label{4.6 log comp}
    \log d_{\widetilde{L}}\le 2\log d_{L} +[L:K]\left(n_{K}\log 4+ \log |N_{K/\mathbb{Q}}(\alpha)|\right).
\end{equation}
If $\sqrt{\alpha}\in L,$ then $\widetilde{L}=L.$ Since $d_{L}$ and $|N_{K/\mathbb{Q}}(\alpha)|$ are positive integers, we have $$d_{\widetilde{L}}=d_{L}\le d_{L}^{2}4^{[L:K]n_{K}}|N_{K/\mathbb{Q}}(\alpha)|^{[L:K]}.$$ So equation \eqref{4.6 log comp} remains valid. Now by using equations \eqref{4.6 log comp}, \eqref{4.4 log disc}, and \eqref{4.2 log norm} we obtain 
\begin{align*}
    \log d_{\widetilde{L}}\le &  384 \log d_{K} +382 \log N_{K/\mathbb{Q}}\mathfrak{m}\\
    &+ 384 n_{K}\log 384 +384 n_{K}\log n_{K}\\
    &+96 n_{K}\log d_{K}+384n_{K}\log N_{K/\mathbb{Q}}\mathfrak{m}\\
    ={}&
    (382+384n_{K})\log N_{K/\mathbb{Q}}\mathfrak{m}+(384+96n_{K})\log d_{K}+384n_{K}\log(384n_{K})\\
    :=&{} R
\end{align*}
Set $T:=1.336 R+1324.8 n_{K}+45.7,$ since $n_{K}\ge1, d_{K}\ge1,$ and $N_{K/\mathbb{Q}}\mathfrak{m}\ge 1,$ we have $$ R\ge 384 n_{K}\log (384n_{K})>384n_{K}\log 384>1920n_{K}.$$ Therefore, $T>3889.92 n_{K}+45.7.$ In particular, $3889.92 n_{K}+45.7$ is larger than both $405+960 n_{K}$ and $1802.3+725.76 n_{K}.$ Moreover, $\log d_{\widetilde{L}}\le R$ gives $$ 1.336\log d_{\widetilde{L}}+1324.8 n_{K}+45.7\le T.$$ So each of the bounds in equation (\ref{4.3 norm prime ideal}) is at most $T^{2}.$ Therefore, $N_{K/\mathbb{Q}}\mathfrak{p}\le T^{2}.$ Now expand $T,$ 

\begin{align*}
    T=&1.336\bigl((382+384n_{K})\log N_{K/\mathbb{Q}}\mathfrak{m}+(384+96n_{K})\log d_{K}\\
    &+384 n_{K}\log (384n_{K})\bigr)+1324.8 n_{K}+45.7\\
    =&\underbrace{\bigl(510.352+513.024n_{K}\bigr)}_{:=C_{1}(K)}\log N_{K/\mathbb{Q}}\mathfrak{m}\\
    &+\\
&\underbrace{
\bigl(513.024+128.256n_{K}\bigr)\log d_K
+513.024n_{K}\log(384n_{K})
+1324.8n_{K}+45.7}_{:=C_2(K)}.
\end{align*}
Consequently, $$ N_{K/\mathbb{Q}}\mathfrak{p}\le \bigl(C_{1}(K)\log N_{K/\mathbb{Q}}\mathfrak{m}+C_{2}(K)\bigr)^{2}.$$
\end{proof}

\section{Proof of the Main Theorem}
In this section, we finally prove the main theorem. We first establish some properties of quadratic characters.
\subsection{Quadratic characters}
Let $\chi_D$ denote the quadratic Galois character, that is
\begin{align*}
&\chi_D:G_K\rightarrow\{\pm 1\}\\
    \chi_D(\sigma)&=
    \begin{cases}
    1, &\quad \text{if }\; \sigma(\sqrt{D})=\sqrt{D},\\
    -1, &\quad \text{if }\; \sigma(\sqrt{D})=-\sqrt{D}.
    \end{cases}    
\end{align*}
\begin{lem}
    $\chi_D$ is unramified at $\mathfrak{p}$ if and only if $K(\sqrt{D})/K$ is unramified at $\mathfrak{p},$ where $\mathfrak{p}$ is a prime ideal of $\mathcal{O}_{K(\sqrt{D})}$ lying over $p$.
\end{lem}
\begin{proof}
    Suppose that $\chi_D$ is unramified at $\mathfrak{p},$ that is, $\chi_D(I_{\mathfrak{p}})=1$. We also have $\text{res}_D(I_{\mathfrak{p}})=I_{\mathfrak{p}}(D)$, and hence, $\# I_{\mathfrak{p}}(D)=1.$
    This implies that $e(K(\sqrt{D}/K)=1$ since $e(K(\sqrt{D}/K)=\# I_{\mathfrak{p}}(D).$ Hence, $\chi_D$ is unramified at $\mathfrak{P}.$\\
    Conversely, assume that $K(\sqrt{D})/K$ is unramified at the prime $\mathfrak{p}.$ Given $\sigma\in I_{\mathfrak{p}}$, let $\overline{\sigma}$ denote its restriction to $I_{\mathfrak{p}}(D).$ Since $I_{\mathfrak{p}}(D)=1$, we have $\overline{\sigma}(\sqrt{D})=\sqrt{D}$, which implies that $\chi_D(\sigma)=1$. Therefore, $\chi_D$ is unramified at $\mathfrak{p}.$
    
\end{proof}

If $\mathfrak{p}\nmid D$, then 
\begin{equation}\label{genlegendre}
    \left(\dfrac{\mathfrak{p}}{K(\sqrt{D})/K}\right)(\sqrt{D})=\left(\dfrac{D}{\mathfrak{p}}\right)_2\sqrt{D}, 
\end{equation}
where the term on the left is the Artin Symbol and on the right we consider the generalized Legendre symbol. (See  \cite[Exercise 5.13]{davidcox} for precise definitions.)
\subsection{Quadratic Twists}
Now consider an elliptic curve $E/K$ given by the Weierstrass equation
\begin{equation*}
    y^2+a_1xy+a_3y=x^3+a_2x^2+a_4x+a_6
\end{equation*}
and let $E_D$ be the quadratic twist of E by $\chi_D$ (See \cite[p.48-49]{COMALADA199445} for more on quadratic twists.) The discriminant of $E_D$ is $D^6$ times the discriminant of $E$. From this we can conclude that 
\begin{equation}\label{rad divides rad}
    \mathrm{rad}(\mathfrak N_{E_D}) \quad \text{divides} \quad \mathrm{rad}(2\mathfrak{N}_ED).
\end{equation}
If we consider a prime $p\nmid N_{K/\mathbb{Q}}(D\mathfrak{N}_E)$, and let $\mathfrak{p}$ be any prime lying over $p$, we have by \cite[Exercise 4.10]{wash2008} that 
\begin{equation*}
    a_{\mathfrak{p}}(E_D)=\left(\dfrac{D}{\mathbb{F}_q}\right)_2a_\mathfrak{p}(E),
\end{equation*}
where $q=p^f$ where $f$ denotes the inertia degree of $\mathfrak p$.
From equation \eqref{genlegendre} above, we have,
\begin{equation*}
    a_{\mathfrak{p}}(E_D)=\chi_D(\text{Frob}_\mathfrak{p})a_\mathfrak{p}(E).
\end{equation*}
\begin{lem}\label{E and ED}
    If $D\neq 1$ is squarefree and $E/K$ has no complex multiplication, then $E$ and $E_D$ are not $K$-isogenous.
\end{lem}
\begin{proof}
    The density of primes $\mathfrak{p}\in\mathcal{O}_K$ for which $\chi(Frob_{\mathfrak{p}})=-1$ is $\frac{1}{2}$ by the Chebotarev density theorem. For any such prime, we either have $a_\mathfrak{p}(E)=-a_\mathfrak{p}(E_D)$ or $a_\mathfrak{p}(E)=0.$ By a generalization of a theorem by Serre, we see that the density of primes for which $\mathfrak{p}\nmid 2\mathfrak{N}_ED$ and $a_\mathfrak{p}(E)=0$ is 0 \cite[Corollary 6.2]{charles2004complexmultiplicationtestselliptic}. Thus, there exists a prime $\mathfrak{p}\nmid 2\mathfrak{N}_ED$ such that $a_\mathfrak{p}(E)=-a_\mathfrak{p}(E_D)$ and hence $E$ and $E_D$ are not $K$-isogenous.
\end{proof}
\subsection{A Special character}  The  construction of this character dates back to Serre \cite[p. 295]{Serre1971}. We consider the following Galois character given by the composition
\begin{equation}\label{quadratic Galois character}
    \varepsilon_\ell:\gal\rightarrow \text{Im}(\rho_{E,\ell})\rightarrow \dfrac{C_{ns}^+(\ell)}{C_{ns}(\ell)}\xrightarrow{\sim}\{\pm 1\}.
\end{equation}
\\
\begin{lem}\label{character properties}
Let $\varepsilon_\ell$ be as above. We have the following properties
     \begin{enumerate}[label={(\roman*)},itemindent=1em]
            \item For each prime $\mathfrak{p}\nmid \mathfrak{N}_E$, $\varepsilon_\ell$ is unramified at $\mathfrak{p}.$
            \item One has $\varepsilon_\ell=\chi_D$ for some $D\mid \mathfrak{N}_E.$
            \item For each prime $\mathfrak{p}\nmid \mathfrak{N}_E$, if $\varepsilon_\ell(\text{Frob}_\mathfrak{p})=-1$, then $a_\mathfrak{p}(E)\equiv 0 \pmod  \ell.$
     \end{enumerate}
\end{lem}
\begin{proof}
We refer the reader to \cite[Lemma 15, p. 1412]{maylewang} for a proof. In our case, we replace the prime $p$ with a prime ideal $\mathfrak p$ and the conductor $N_E$ with the conductor ideal $\mathfrak N_E$. By $D\mid\mathfrak N_E$ we mean a containment of ideals, that is, $\mathfrak N_E \subset (D) $.
\end{proof}
Now that we have established the basic facts, we prove the main theorem
\subsection{Proof of the Main Theorem}
\begin{theorem}\label{5.4mainthm}
    Assume GRH. Let $E$ be an elliptic curve without complex multiplication and conductor ideal $\mathfrak N_E$, over a Galois number field $K$. For a prime number $\ell$ with
    \begin{equation*}
        \ell\geq 2C_1(K) \log N_{K/\mathbb Q}\operatorname{rad}\left( 2\mathcal O_K\mathfrak N_{E}\right )+2C_2(K). 
    \end{equation*}
    $\rho_{E,\ell}(G_K)$ is not contained in the normalizer of a non-split Cartan subgroup of $\mathrm{GL}_2(\mathbb{F}_\ell)$.
\end{theorem}
\begin{proof}
     Assume that $G_E(\ell)$ is contained in the normalizer of a non-split Cartan subgroup $C_{ns}^+(\ell)$ but not fully contained in said non-split Cartan subgroup $C_{ns}(\ell)$ (up to conjugation) Consider the quadratic Galois character given in equation \eqref{quadratic Galois character}. Let $D$ be a square free algebraic integer such that $\varepsilon_\ell=\chi_D$ and consider the quadratic twist $E_D$ of E. By Lemma~\ref{character properties}(ii) and (\ref{rad divides rad}) we have that
    \begin{equation}
        \mathrm{rad}(\mathfrak N_{E_D})\quad\text{divides}\quad \mathrm{rad}(2\mathfrak{N}_E).
    \end{equation}
    For each prime ideal $\mathfrak{p}\nmid 2\mathfrak{N}_E$, we have that 
    \begin{equation}\label{relation between traces}
        a_\mathfrak{p}(E)=\varepsilon_\ell(\text{Frob}_\mathfrak{p})a_\mathfrak{p}(E_D).
    \end{equation}
    By Lemma \ref{E and ED}, $E$ and $E_D$ are not K-isogenous. Thus, there exists a prime ideal $p\nmid 2\mathfrak{N}_E$ such that 
    \begin{equation}\label{trace unequal E ED}
        a_\mathfrak{p}(E)\neq a_\mathfrak{p}(E_D).
    \end{equation}
    Let $\mathfrak{p}_0$ be the prime ideal of least norm such that $\mathfrak{p}_0\nmid 2\mathfrak{N}_E$ and $(\ref{trace unequal E ED})$ holds. Applying Theorem \ref{Eff iso for K} to $E$ and $E_D$ and using equation (\ref{rad divides rad}), we find that 
    \begin{equation}
       N_{K/\mathbb Q}\mathfrak p_0\le\left(C_1(K) \log N_{K/\mathbb Q}\operatorname{rad}\left( 2\mathcal O_K\mathfrak N_{E_1}\right )+C_2(K) \right)^2,
    \end{equation}
    where the constants are explicitly described in Theorem \ref{Eff iso for K}.
    Considering (\ref{relation between traces}) and (\ref{trace unequal E ED}), we see that $\varepsilon_\ell(\text{Frob}_{\mathfrak{p}_0})=-1.$ Thus by lemma \ref{character properties}(iii), $\ell$ divides $a_{\mathfrak{p}_0}(E),$ and so $\ell\leq |a_{\mathfrak{p}_0}(E)|.$ The Hasse bound gives us that $|a_{\mathfrak{p}_0}(E)|\leq 2\sqrt{N(\mathfrak{p}_0})$. Thus,
    \begin{equation*}
        \ell\leq 2\sqrt{N(\mathfrak p_0)}\leq 2\left(C_1(K) \log N_{K/\mathbb Q}\operatorname{rad}\left( 2\mathcal O_K\mathfrak N_{E_1}\right )+C_2(K) \right).
    \end{equation*}
    Expanding the RHS gives us the claimed bound.
\end{proof}

\bibliographystyle{plain}
\bibliography{sample.bib}
\end{document}